\documentstyle[xypic,12pt]{article}
\newtheorem{th}{\hspace{\parindent}Theorem}[section]
\newtheorem{lem}[th]{\hspace{\parindent}Lemma}
\newtheorem{prop}[th]{\hspace{\parindent}Proposition}

\newtheorem{cor}[th]{\hspace{\parindent}Corollary}

\newtheorem{ass}[th]{\hspace{\parindent}Assertion}

    \newcommand{\dfi}{\par\refstepcounter{th}%
    {\bf Definition \thesection.\arabic{th}\quad}}
    \newcommand{\rem}{\par\refstepcounter{th}%
    {\bf Remark \thesection.\arabic{th}\quad}}
    
    \newcommand{\cst}{\par\refstepcounter{th}%
    {\bf Construction \thesection.\arabic{th}\quad}}
    \newcommand{\pri}{\par\refstepcounter{th}%
    {\bf Example \thesection.\arabic{th}.\quad}}
\def\Dok{{\bf Proof: \/}}

\def\o{\otimes}
\def\La{\Lambda}
\def\la{\lambda}
\def\f{\varphi}
\def\vfi{\varphi}
\def\w{\wedge}

\def\p{\psi}
\def\f{\varphi}
\def\io{\iota}

\def\te{\theta}
\def\t{\tau}
\def\ro{\rho}
\def\e{\varepsilon}
\def\ss{\subset}
\def\s{\sigma}
\def\x{\chi}

\def\r{\rho}
\def\bop{\oplus}
\def\Lra{\Leftrightarrow}
\def\hop{\widehat{{}\oplus {}}}
\def\b{\beta}
\def\s{\sigma}
\def\a{\alpha}
\def\g{\gamma}
\def\G{\Gamma}
\def\F{\Phi}

\def\<{\langle}
\def\>{\rangle}

\newfont{\pz}{msbm10 at 12pt}
\def\skr{\mbox{\pz o}\,}

\newfont{\gt}{eufm10 at 12pt}
\def\gg{{\mbox{\gt g}}}

\def\P{{\cal P}}
\def\M{{\cal M}}
\def\N{{\cal N}}

\def\H{{\cal H}}
\def\K{{\cal K}}
\def\A{{\cal A}}
\def\B{{\cal B}}
\def\E{{\cal E}}

\def\c{{\scriptscriptstyle{\circ\atop\vphantom{a}}}}
    \newcommand{\operatorname}[1]{\,{\rm #1}\,}

\def\Aut{\operatorname{Aut}}

\def\ind{\operatorname{ind}}
\def\index{\operatorname{index}}

\def\End{\operatorname{End}}

\def\Tor{\operatorname{Tor}}

\def\Vect{\operatorname{Vect}}

\def\Mor{\operatorname{Mor}}
\def\OP{\operatorname{OP}}
\def\Hom{\operatorname{Hom}}

\def\dim{\operatorname{dim}}

\def\supp{\operatorname{supp}}

\def\CZ{\operatorname{CZ}}
\def\Int{\operatorname{Int}}
\def\Smbl{\operatorname{Smbl}}

\def\Coker{\operatorname{Coker}}
\def\Ker{\operatorname{Ker}}

\def\Tor{\operatorname{Tor}}
\def\Proj{\operatorname{Proj}}

\def\Ob{\operatorname{Ob}}
\def\Id{\operatorname{Id}}
\def\Im{\operatorname{Im}}

\def\ai{\operatorname{a\mbox{-}ind}}
\def\ti{\operatorname{t\mbox{-}ind}}
\def\ci{\operatorname{{\scriptscriptstyle \C}\mbox{-}ind}}
\def\R{{\bf R}}

\def\Z{{\bf Z}}
\def\C{{\bf C}}

\def\Lra{\Leftrightarrow}
\def\lra{\longrightarrow}
\def\wt{\widetilde}
\def\ola#1{\stackrel{#1}{\longrightarrow}}

\def\hop{\widehat{\oplus}}
\def\b{\beta}
\def\a{\alpha}
\def\g{\gamma}

\def\P{{\cal P}}

\def\Ker{\operatorname{Ker}}
\def\Id{\operatorname{Id}}
\def\Im{\operatorname{Im}}
\def\GL{\operatorname{GL}}

\def\wA{{\widetilde\A}}
\def\wB{{\widetilde\B}}

\def\cC{{\cal C}}
\def\bL{{\bf\La}}

\def\bK{{\bf K}}
\def\bVect{{\bf Vect}}
\def\bL{{\bf L}}
\def\bM{{\bf M}}

\def\ot{\otimes}
\def\otl{\mathop{\otimes}\limits}

\def\wo{\widehat{\otimes}}

\newcommand{\Mat}[4]{\left( \begin{array}{cc}
                            #1 & #2 \\
                            #3 & #4
                      \end{array} \right)}

\author{
E.~V.~Troitsky\thanks
{The work is partially supported by the RFBR grant
N 99-01-01201
and INTAS grant N 96-1099.
}
}
\date{}
\title{
Actions of compact groups,
$C^*$-index theorem, and families
}
\begin{document}
\maketitle
\begin{center}
Moscow State University\\
Faculty of Mechanics and Mathematics \\
Chair of Higher Geometry and Topology\\
119 899 Moscow  Russia\\
troitsky@mech.math.msu.su\\
Home page URL: http://mech.math.msu.su/\~{}troitsky
\end{center}

\begin{abstract}
We prove the index theorem for elliptic operators acting on sections of
bundles where fiber is equal to a projective module over a $C^*$-algebra,
in the situation of action of a compact Lie group on this algebra as well
as on the total space commuting with symbol.  As an application the
equivariant index theorem for a direct product of base by the space of
parameters is obtained.

The present preprint is a detailed version of a combination of two
papers to appear in "J. Math. Sci." and "Ann. Global Anal. Geom."
\end{abstract}

\tableofcontents

 \section{Introduction}
 \markright{}

The rational version of index theorem for $C^*$-elliptic operators
\cite{Mis,MF} has numerous applications in differential topology (to
proof the Novikov conjecture), differential geometry (curvature of
spin manifolds), theory of elliptic operators with random and
periodic coefficients, etc. After that the equivariant generalization
of this theorem \cite{TroIAN}, as well as equivariant generalization
taking the torsion into account \cite{TroStekl} was obtained.  They
also have found their applications (see, e.g. \cite{RosWein,SolTro}).
In these theorems a compact Lie group acts on the total space of a
bundle of $A$-modules, but not on the algebra $A$ itself. For some
applications, the most classic of which is calculation of equivariant
index of a family of elliptic operators, it is necessary for the
group to act on $A$ compatible with the action on the total space.
The present paper is devoted to the decision of this problem.

We use intensively the foundations of theory of
$C^*$-Hilbert modules. One can find it in~\cite{Lance,WO,MaTroJMS}.
The modern state of index theory for families and its applications
can be found e.g. in \cite{Tsuboi,Bism}.

\smallskip
The author is grateful to V.~M.~Manuilov and
A.~S.~Mishchenko for helpful discussions.

 \section{An averaging theorem}

Let $G$ be a compact Lie group acting  continuously on $C^*$-algebra
$A$ by involutive automorphisms. If $A$ is unital, then the unity has
to be invariant. Then $g\in G$ takes self-adjoint elements to
self-adjoint and positive ones to positive.

\rem\label{rem:posit}
If $a\ge b\ge 0$, then $ga\ge gb$. Indeed, $a-b=c^*c\ge 0$. Hence,
$ga-gb=g(c)^*g(c)\ge 0$.

An $A$-module $M$ is called $GGA$-{\em module\/} if it is equipped
with a $\C$-linear action of $G$, such that
 $$
g(m\cdot a)=g(m)\cdot g(a),\qquad m\in M,\qquad a\in A.
 $$

\dfi\label{def:invproizv} An inner $A$-valued product on Hilbert
module is called {\em invariant} or $GGA$-product if
 $$
\<gx,gy\>=g(\<x,y\>).
 $$

Let us remark that the averaging theorem of \cite{TroIAN} (see also
\cite{SolTro}) does not define an invariant product in the $GGA$-case.
Moreover, the formula
$$ \<u,v\>\,=\,\int _G \,\<T_xu,T_xv\>^\prime\,dx, $$
does not define an $A$-inner product:

\pri Let $\M=A=C(G)$. Then the new ``product'' is valued in constants
only. Some effects of such originality are studied in \cite{FMTZAA}.

\dfi Let us define an action of $G$ on the module $l_2(A)$ over a
unital $G$-algebra $A$ by the formula
$g(u_1,u_2,\dots)=(gu_1,gu_2,\dots)$.

 \begin{lem}\label{lem:spezdejstvije}
This is a continuous action, and the initial $A$-valued product is
invariant with the respect to it.
 \end{lem}

\Dok
First of all
\begin{eqnarray*}
g(ua)&= &
g((u_1,u_2,\dots)a)=(g(u_1a),g(u_2a),\dots)
=(g(u_1)g(a),g(u_2)g(a),\dots)=
\\
&=&(g(u_1),g(u_2),\dots)g(a)=g(u)g(a).
\end{eqnarray*}
Further,
 $$
\<gu,gv\>=\sum_{i=1}^\infty (g(u_i))^* g(v_i)=
\sum_{i=1}^\infty g(u^*_i) g(v_i)=
\sum_{i=1}^\infty g(u^*_iv_i)=
g\left(\sum_{i=1}^\infty u^*_iv_i\right)=g(\<u,v\>).
 $$
Let us demonstrate the continuity. Let $\|g-h\|<\e$ as automorphisms
of the algebra $A$. Then
\begin{eqnarray*}
\|\<(g-h)u,(g-h)u\>\|&=&
\left\|\sum_{i=1}^\infty ((g-h)(u_i))^* (g-h)(u_i)\right\|
=\left\|\sum_{i=1}^\infty (g-h)(u_i^* u_i)\right\|=\\
&=&\left\|(g-h)\sum_{i=1}^\infty u_i^* u_i\right\|\le \e \|u\|^2.
\qquad \Box
\end{eqnarray*}

     Let us   recall   some facts about the integration  of
operator-valued functions \cite{Kas1}. Let  $X$  be   a
compact space, $A$  be a $C^*$-algebra,  $\vfi :  C(X) \to A$ be an
involutive homomorphism of algebras with unity, and $F:X\to A$ be
a continuous map, such that for every $x\in X$ the element $F(x)$
commutes with the image  of  $\vfi$.  In  this  case  the
integral
      $$
     \int_X F(x)\,d\vfi \quad \in \quad A
      $$
can be defined in the following way.
Let $X=\cup_{i=1}^n U_i $ be an open covering and
$$
\sum_{i=1}^n\a _i(x)=1
$$
be a corresponding partition of  unity. Let us choose the points
     $\xi _i \in U_i$ and compose the integral sum
 $$
\sum (F, \{U_i\}, \{\a _i\}, \{\xi _i\}) = \sum _{i=1}^n F(\xi _i)
     \vfi (\a _i).
      $$
     If there  is  a  limit  of  such integral sums then it is
called the corresponding integral.

If $X$ is a Lie group $G$ then it is natural to take $\f$ equal to
the Haar measure
      $$
     \f : C(X) \to \C, \quad \f (\a ) = \int _G \a (g)\,dg
      $$
(though this is only a positive linear map, not a $*$-homomorphism)
     and to define for a norm-continuous $Q: G \to L(H)$
($A$ is realized as a subalgebra in algebra $L(H)$ of all bounded
operators on Hilbert space $H$)
      $$
\int _G  Q(g)\, dg = \lim \sum _i Q(\xi _i) \int _G \a _i (g)\,
dg.
      $$

   If $Q: G\to P^ {+} (A) \ss L (H) $, then, since
      $
     \int _G \a _i (g) \, dg \ge 0,
      $
we get
      $$
     \sum _i Q (\xi _i) \cdot \int _G \a _i (g) \, dg \; \in
\; P^ {+} (A)
\quad {\rm and} \quad
    \int _G Q (g) \, dg \; \in \; P^ {+} (A)
     $$
 (the positive cone $ P^ {+} (A) $ is convex and closed).
So we have proved
the following Lemma.

\begin{lem}\label{lem:aver46}
Let $Q: G\to P^{+}(A)$ be a continuous function.
Then for the integral in  the  sense  of
{\rm \cite{Kas1}} we have
      $$
     \int _G Q (g) \, dg \ge 0. \qquad \Box
      $$
\end{lem}

 \begin{th}\label{usred}
  Let $\GL=\GL (A)$ be
\label{ind:gl}
the full general linear group, i.~e.,
  the group of all bounded $A$-linear automorphisms of $l_2(A)$.
  Let $g \mapsto T_g \quad (g\in  G,  T_g  \in  \GL )$  be  a
  representation of $G$ such that the map
      $$
     G \times l_2(A)\to l_2(A), \quad (g,u)\mapsto T_gu
      $$
  is continuous.
   Then on  $l_2(A)$  there is a $GGA$-product equivalent to
  the original one.
     \end{th}

\Dok
Let $\< \,,\,\>^\prime$ be the original product.
      We have  a  continuous   map
$ G\to   A,\quad   x\mapsto x^{-1}(\<T_x u,T_x v\>^\prime).$
 for every $u$ and $v$ from $l_2(A)$.  We define
 the new product by
      $$
     \<u,v\>\,=\,\int _G \,x^{-1}(\<T_x u,T_x v\>^\prime)\,dx,
      $$
     where the  integral can be defined in the sense of either of the two
     definitions from \cite[p. 810]{Kas1}   because   the  map  is
     continuous with the respect to the norm of the $C^*$-algebra.
This product is an $A$-Hermitian map
$l_2(A)\times l_2(A)\to  A$. Indeed,
\begin{eqnarray*}
\<u\cdot a,v\cdot b\>&=&
\,\int_G \,x^{-1}(\<T_x (u\cdot a),T_x (v\cdot b)\>')\,dx=\\
     &=&
\,\int_G \,x^{-1}(\<T_x (u)\cdot (xa),T_x (v)\cdot (xb)\>')\,dx=\\
&=&
\,\int_G \,x^{-1}\Bigl[(xa)^*(\<T_x (u),T_x (v)\>')(xb)\Bigr]\,dx=\\
     &=&
\,\int_G \,x^{-1}\Bigl[x(a^*)(\<T_x (u),T_x (v)\>')(xb)\Bigr]\,dx=\\
&=&
\,\int_G \,(a^*)x^{-1}(\<T_x (u),T_x (v)\>')b\,dx=
a^*\<u,v\>b.
\end{eqnarray*}
Since $x^{-1}:A\to A$ is an involutive mapping, it takes positive
elements to positive ones and by Lemma \ref{lem:aver46}
 $$
\<u,u\>=\,\int _G \,x^{-1}(\<T_x u,T_x u\>^\prime)\,dx,\ge 0.
 $$
For
$
T_x u=(a_1(x),a_2(x),\dots) \in l_2(A)
$
the equality $\<u,u\>=0$ takes the form
      $$
     \int _G \sum _{i=1}^\infty x^{-1}(a_i(x) a_i^*(x))\, dx =0.
      $$
Since $x^{-1}(a_i(x) a_i^*(x))\ge 0$, we have $a_i(x)=0$ almost
everywhere. Hence, $a_i(x)=0$ for all
$x$ by the continuity, and $T_x u=0$. In particular,
$u=0$.

Let us demonstrate the continuity of this new product. Fix
$u\in l_2(A)$. Then by Lemma \ref{lem:spezdejstvije},
      $     x\mapsto T_x(u),\quad G \to l_2(A)      $
is a continuous mapping of a compact space. Hence, the set
$\{T_x(u)\,|\, x\in G  \}$ is bounded. So, by the principle
     of uniform boundness  \cite[vol.~2, p.~309]{BaRo}
\begin{equation}\label{aver16}
     \lim _{v\to 0} T_x (v) = 0
\end{equation}
uniformly with respect to $x\in G$. If $u$ is fixed then
      $$
     \| T_x(u)\| \le M_u = const
      $$
and by (\ref{aver16})
\begin{eqnarray*}
\| \<u,v\>\|&=&\|\int_G \,x^{-1}(\<T_x u,T_x v\>^\prime)\,dx\| \le \\
&\le& M_u \cdot vol\,G\cdot \sup\limits_{x\in G} \| T_x (v)\| \to 0
\quad (v\to 0).
 \end{eqnarray*}
This gives the continuity at  0 and hence
on the whole space $l_2 (A) \times l_2 (A) $. Since each operator
$T_y$ is an automorphism, we get
\begin{eqnarray*}
\<T_y u,T_y v\>&=&\int _G \,x^{-1}(\<T_x T_y u,T_x T_y v\>^\prime)\,
dx=\\
&=&\int _G \,y(xy)^{-1}(\<T_{xy} u,T_{xy} v\>^\prime)\,dx=\\
&=&y(\<u,v\>),
\end{eqnarray*}
i.~e. this product is invariant.

Now we  will  show  the equivalence of the two norms and,
     in particular,  the continuity of  the  representation.  There
     is a number $N>0$ such that $\| T_x\|^\prime \leq N$ for every
     $x\in G$. So by \cite{Kas1}
(for the simplicity $vol\:G=1$)
\begin{eqnarray*}
\| u\|^2 & =&\|\<u,u\>\|_A
     =\|\int_G\,x^{-1}(\<T_x u,T_x u\>^\prime)\,dx\|_A\le \\
     &\le &\left( \sup _{x\in G}
     \| T_x u\|^\prime\right) ^2 \le N^2 (\| u\|^\prime ) ^2.
\end{eqnarray*}
On the other hand, let
$\<u,u\>'=1$. Then applying Lemma \ref{lem:aver46} and
Remark \ref{rem:posit}, we obtain
      \begin{eqnarray*}
\<u,u\>^\prime &=&\int _G g^{-1}(\<u,u\>)\, dg=
     \int _G g^{-1}(\<T_{g^{-1}}
       T_g  u,T_{g^{-1}} T_g u\>^\prime)\, dg \le \\
&\le& \int _G g^{-1}(\| T_{g^{-1}}\|^2 \<T_g  u, T_g u\>')\, dg \le \\
      &\le& \int _G g^{-1}(N^2\<T_g  u_g,
       T_g u_g\>^\prime)\, dg
     = N^2 \int _G g^{-1}(\<T_g  u, T_g u\>^\prime)\, dg =
      N^2 \<u,u\>.
      \end{eqnarray*}
Then
      $
      (\| u\|^\prime)^2  =
     \|\<u,u\>^\prime\|_A  \le N^2 \|\<u,u\> \|_A = N^2\| u\| ^2.
     $
By linearity we obtain a similar estimate for $u$ with invertible
$\<u,u\>'$, while the elements of such a form are dense
(see \cite{DupFil}) in $l_2(A)$.  \hfill$\Box$

\rem\label{usrP}
    $l_2(P)$ is a direct summand in $l_2(A)$, so
the previous theorem holds for
     $l_2(P)$.


    \section{$K$-theory of $GGA$-bundles }\label{subs:KA}

Let us recall some general constructions of K-theory, contained
in \cite{KaKn}, and also \cite{KarCli}.

\dfi\label{kar1}~\cite[I.6.7]{KaKn}
An additive category $\cC$ is called
{\em pseudo-Abelian\/},
\label{ind:catpseud}
 if for each object $E$ from $\cC$ and each
morphism $p:E\to E$, satisfying to a condition $p^2=p$ (i.~e. an
idempotent) there exists the kernel $\Ker p$.
For an arbitrary additive
category $\cC$ there exists
{\em associated pseudo-Abelian category\/}%
\label{ind:catpseudass}
$\wt\cC$ which is a solution of the appropriate
universal problem and is defined as follows \cite[I.6.10]{KaKn}.
Objects of $\cC$ are pairs $(E, p)$, where $E\in\Ob (\cC) $
and $p$ is a projector in $E$. A morphism from $(E, p)$ to $(F, q)$
is such a morphism $f: E\to F$ of the category $\cC$, that
$f\circ p=q\circ f=f$.

\dfi\label{kar2}~\cite[\S~II.1]{KaKn}
We call {\em simmetrization\/}
\label{ind:simmetrization}
  of an Abelian monoid $M$ the
following Abelian group $S (M)$.
Consider the product $M\times M$ and its
quotient monoid with respect to
the equivalence relation
 $$
(m,n)\sim (m',n')\:\Lra\:\exists\:p,q:\:(m,n)+(p,p)=(m',n')+(q,q).
 $$
This quotient monoid is a group denoted $S(M)$.
If we consider now an
additive category $\cC$ and denote through $\dot E$ the
isomorphism  class of an
object $E$ from $\cC$, then the set $\F (\cC) $ of these classes
is equipped with
a structure of an Abelian monoid with respect to
operation $\dot E + \dot F= (E\oplus F)^\bullet$.
In this case the group $S (\F (\cC)) $ is denoted
through $K (\cC) $ and is called
{\em Grothendieck group\/}
\label{ind:GrothGroup}
of the category $\cC$.

\dfi\label{kar3}~\cite[\S~II.2] {KaKn}
{\em Banach structure\/} on an additive category $\cC$ is
defined by the introducing of the structure
of a Banach space on all groups $\cC (E, F) $, where $E$ and $F$ are
arbitrary objects from $\cC$. It is assumed, that
applications of
composition of morphisms $\cC (E, F) \times\cC (F, G) \ola {} \cC (E, G)
 $
are bilinear and continuous. In this case we call $\cC$ {\em Banach
category\/}.
\label{ind:catbanach}

\dfi\label{kar4}~\cite[\S~II.2]{KaKn}
Let $\cC$ and $\cC'$ be additive categories. An additive functor
$\f:\cC\to\cC'$ is called {\em quasi-surjective\/}
\label{ind:functquasi}
 if each object
of $\cC'$ is a direct summand of an object of type $\f (E)$.
A functor
 $\f$ called {\em full\/}
\label{ind:functpoln}
if for any $E,\,F\in\Ob (\cC)$
the map $\f (E, F) :\cC (E, F) \to\cC' (\f (E),\f (F)) $ is
surjective. For Banach categories $\f$ is called
{\em Banach\/}
\label{ind:functbanach}
if this map $\f (E, F)$ linear and continuous.

\dfi\label{kar5}~\cite[II.2.13]{KaKn}
Let $\f:\cC\to\cC'$ be a quasi-surjective Banach functor. We shall
denote by $\G(\f)$ the set consisting of triples of the form
$(E, F,\a)$, where $E$ and $F$ are objects of the category $\cC$ and
$\a:\f (E) \to\f (F) $ is an isomorphism.  The triples $ (E, F,\a) $
and $(E',F',\a') $ are named {\em isomorphic\/}, if there are such
isomorphisms $f: E\to E'$ and $g: F\to F'$, that the diagram
 $$
\diagram
\f(E)\rto^\a \dto_{\f(f)}& \f(F)\dto^{\f(g)}\\
\f(E')\rto^{\a'}&\f(F')
\enddiagram
 $$
commutes. A triple $(E,F,\a)$ is {\em elementary\/}
if $E=F$ and isomorphism $\a$ is homotopic in the set of
automorphisms of
$\f (E) $ to the identical isomorphism $\Id_{\f (E)}$.
We define
{\em sum\/} of two triples $ (E, F,\a) $ and $ (E ', F ',\a ')$
as
$$
(E\oplus E ', F\oplus F ',\a\oplus \a ').
$$
{\em The Grothendieck group $K (\f)$ of a functor $\f$\/}
is defined as quotient set of the monoid $\G(\f)$
with respect to the following
equivalence relation:  $\s\sim\s'$ if and only if there exist such
elementary triples $\t$ and $\t'$, that the triple $\s + \t$ is
isomorphic to the triple $\s' + \t'$.  The operation of addition
introduces on $K (\f) $ a
structure of Abelian group. The class of a triple we shall denote
by $d (E, F,\a) $.

\dfi\label{kar6}~\cite[II.3.3]{KaKn}
Consider the set of pairs of the form $(E,\a)$, where
$E$ is an object of the category $\cC$ and $\a$ is an automorphism
of $E$.
Two pairs $(E,\a)$ and $(E',\a')$ are called {\em isomorphic\/}, if
there is such isomorphism $h: E\to E'$ in category $\cC$, that the
diagram
 $$
\diagram
E\rto^h\dto_\a& E'\dto^{\a'}\\
E\rto^h&E'
\enddiagram
 $$
commutes. The direct sum defines the
operation of {\em addition\/} of pairs. A pair $(E,\a)$ is
called {\em elementary\/}, if the automorphism $\a$ is homotopic
to $\Id_E$ in the set of automorphisms of $E$. Abelian group
\cite[II.3.4]{KaKn} $K^{-1}(\cC)$ is defined as a quotient set
(with operation of addition) of the set of pairs $\{(E,\a) \}$
with respect to the
following equivalence relation:  $\s\sim\s'$ if and only if
 there
are such elementary pairs $\t$ and $\t'$, that $\s + \t$ is isomorphic
to $\s' + \t'$.

\dfi\label{kar7}~\cite[II.4.1] {KaKn}
Let $\cC$ be a Banach category and $C^{p, q}$ be the Clifford
algebra.  We shall denote by $\cC^ {p, q}$ the category, objects
of which are pairs $(E,\r)$, where $E\in\Ob (\cC) $ and
$\r: C^ {p, q} \to\End (E) $ is a homomorphism of algebras. A morphism
from a pair $ (E,\r) $ to a pair $ (E ',\r') $ is such a
$\cC$-morphism $f:E\to E'$, that $f\circ\r (\la)=\r(\la)\circ f$ for each
element $\la\in\cC^{p, q} $.

\dfi\label{kar8}~\cite[III.4.11]{KaKn}
Let $\cC$ be a pseudo-Abelian Banach category.
{\em The group\/} $K^{p,q}(\cC)$ is defined as the
Grothendieck group of the functor
$\cC^ {p, q + 1} \to\cC^ {p, q} $ in the sense
of the Definition~\ref{kar5}.

The following statement can be easily obtained by
the properties of Clifford algebras.

\begin{th}\label{kar9}~{\rm \cite[III.4.6, III.4.12]{KaKn}}
The groups $K^{p,q}(\cC)$ depend only on the difference $p-q$. Besides,
the groups $K^{0,0} (\cC)$ and $K^{0,1} (\cC) $ are canonically
isomorphic to groups $K (\cC)$ and $K^{-1}(\cC)$.
\end{th}

\dfi\label{kar9a}
Now we can define
$K^{p-q} (\cC) =K^{p,q} (\cC)$ and similarly for K-groups of functors.

We need also another description of K-groups, which
is equivalent~\cite[\S\S~III.4, III.5]{KaKn} to the initial.

\dfi\label{kar10}~\cite [III.4.11, III.5.1]{KaKn}
Let $\cC$ be a pseudo-Abelian Banach category and
let $E$ be a $C^{p,q}$-module (an object of the category $\cC^{p,q}$).
{\em Gradation\/}
\label{ind:graduirovan}
of $E$ is such endomorphism
$\eta$ of object $E$ (considering as an object from $\cC$), that
\begin{enumerate}
\item $\eta^2=1$,

\item $\eta\r (e_i) =-\r (e_i) \eta$, where $e_i$ are
 the generators of Clifford algebra
and $\r: C^{p,q} \to\End(E)$ is the homomorphism, determining
the $C^{p,q}$-structure on $E$.
 \end{enumerate}

In other words, a gradation of $E$ is a
 $C^{p,q+1}$-structure
on $E$, extending the initial $C^{p,q}$-structure (if we put
$\r(e_{p + q + 1}) =\eta$).

The term ``gradation'' arises from the following fact.
Morphism $\eta$
determines the decomposition into a direct sum
$$
E=E_0\oplus E_1,\quad
E_0=\Ker\left(\frac{1-\eta}2\right), \quad
E_1=\Ker\left(\frac{1+\eta}2\right),
$$
while
$\r: C^ {p, q} \to\End (E_0\oplus E_1) $ is a morphism
of $\Z/2\Z$-graded algebras.

Let us define group $K^{p,q} (\cC) $ as a quotient group
of the free Abelian group,
generated by triples $(E,\eta_1,\eta_2)$,
where $E$ is a $C^{p,q}$-module and $\eta_1,\,\eta_2$ is a gradation
of $E$, with respect to the subgroup, generated by relations
\begin{enumerate}
\item
$(E,\eta_1,\eta_2)\oplus(F,\xi_1,\xi_2)=(E\oplus
F,\eta_1\oplus\xi_1,\eta_2\oplus \xi_2) $,

\item $(E,\eta_1,\eta_2)=0$, if $\eta_1$ is homotopic to $\eta_2$ in
  the set of gradations of $E$.
\end{enumerate}
As usual, by
$d(E,\eta_1,\eta_2) \in K^{p, q} (\cC)$ we shall
denote the class of triple $ (E,\eta_1,\eta_2) $.

We pass to the necessary
specification of these constructions.

If $X$ is a paracompact topological space, let
$\Vect (X; A)$ be the
category of locally trivial bundles $p: E \to X$
with fiber $M=p^{-1} (x) \in \P (A)$ and structure group
equal to $\Aut_A M$.
Such bundles are called $A$-{\em bundles\/}.
This category
is Banach in the sense of Definition~\ref{kar3}
(see also \cite{KarCli,MF,KarTrans}).
Let $G$ be a compact Lie group, acting on $X$ and algebra $A$
continuously.

An $A$-bundle $p:E\to X$ is called
{\em $GGA$-bundle\/}
\label{ind:gabundl}
if a $G$-space structure is given on $E$,
and for any elements $g \in G$, $e \in E$, $x\in X$
\begin{description}
 \item [1.] $gp (e) =pg (e)$, and
 \item [2.] $g: p^ {-1} (x) \to p^ {-1} (gx)$ is an $GGA$-linear
mapping (i.~e. $g(e\cdot a)=g(e)\cdot g(a)$).
\end{description}

We form the Banach category $\bVect_G (X; A)$ whose
objects are $GGA$-bundles and whose morphisms
are the morphisms of $\Vect (X; A) $ that commute with
the action of $G$.

     The set of all continuous sections $s:X \to E$ of a
$GGA$-bundle $E$ over a compact space $X$  forms
a Banach $A$-module $\G(E)$ (in the topology
of the maximum of the norm). The group $G$ strongly
continuously acts on $\G (E) $ according to the rule
$(Gs)(x) =gs(g^{-1} x), $ where $g \in G$, $s \in \G (E) $,
    and $x \in X$.
Let us note that this module is a $GGA$-module over the algebra of all
continuous $A$-valued functions on $X$:
\begin{eqnarray*}
 g(s\cdot a)(x)&=& g((s\cdot a)(g^{-1}x))=g(s(g^{-1}x)\cdot
a(g^{-1}x))=\\
 &=& g(s(g^{-1}x))\cdot  g(a(g^{-1}x)) = g(s)(x)\cdot g(a)(x)=\\
&=&(g(s)\cdot g(a))(x).
\end{eqnarray*}

The averaging mapping $s\mapsto\int_G gs$ determines a projection
$\mu: \G (E) \to \G^G (E)$, where $\G^G (E)$ is the space of
$G$-invariant sections.

\begin{lem}\label{prodoljsec}
Let $s'$ be a $G$-invariant cross-section of a
    $GGA$-bundle $E\to X$ over a closed $G$-stable
subset\footnote{That is $GY\ss Y$.} $Y$ of the compact
    $G$-space $X$. Then $s'$ can be extended
    to a $G$-invariant cross-section over $X$.
\end{lem}

\Dok
Similarly to~\cite[2.1.1]{AK}.
\hfill$\Box$

\smallskip
Consider two $GGA$-bundles $E$ and $F$ over the compact
base $X$.
The $G$-bundle $\Hom(E,F)$ is introduced in the standard way.
It is not an $A$-bundle,
but Lemma~\ref{prodoljsec} remains valid for it.
  As usual, we can identify  $\G^G (\Hom (E, F))$ with the set
  of $GGA$-morphisms $\vfi: E \to F$. Using
Lemma~\ref{prodoljsec}, we
  obtain the following assertion.

\begin{lem}\label{prodoljmor}
Let $\f': E|_Y\to F|_Y$ be a morphism of the
restrictions of $GGA$-bundles $E$ and $F$ over a compact space $X$ to a
closed $G$-stable subset $Y$. Then $\vfi'$ can be extended
to a morphism of $G$-$A$-bundles $\vfi: E \to F$ over
     $X$. If $\vfi'$ is an isomorphism, then
     there exists a $G$-stable open neighborhood $U$ of $Y$
     such that
     $$
     \vfi\vert_U: E\vert_U \to F\vert_U
     $$
is an isomorphism. Any two such extensions
$\vfi_0$ and $\vfi_1$ are homotopic to each other over some
     $G$-neighborhood  $U'\supset Y $ in the class
of isomorphisms.\hfill$\Box$
\end{lem}
If $E$ is a $G$-space and $I$ is the unit interval, then
we define the action of $G$ on $Y \times I$ by the formula
$g (y, t) = (gy, t)$.
  From Lemma~\ref{prodoljmor},
as well as in a classical case, we obtain
the following fact.

\begin{lem}\label{homot}
Let $Y$ be a compact $G$-space,
$f_t: Y \to X$ a homotopy of $G$-mappings $(0\le t\le 1)$,
and $E$ a  $G$-$A$-bundle over $X$. Then $f^*_0 E \cong f^*_1 E$.
   \hfill$\Box$
\end{lem}

By Lemma \ref{prodoljmor} we can define $GGA$-bundle
$$
E_1 \cup _ {\vfi} E_2 \to X
$$
in the usual way, where
$X=X_1 \cup X_2$, $Y=X_1 \cap X_2$, $X_1$ and $X_2$ are closed
$G$-subspaces of the compact space $X$,
  $E_1 \to X_1$ and
$E_2 \to X_2$ are $GGA$-bundles, and
$\vfi: E_1\vert_Y \to E_2 \vert_Y$ is an isomorphism.
Further, as follows from Lemma~\ref{homot},
$E_1\cup_{\vfi} E_2$ up to isomorphism
depends only on the $G$-homotopy class of $\vfi$.

\begin{th}\label{epim} {\rm (see, e.~g.~\cite{SolTro})}
    Let $E$ and $F$ be
 $A$-bundles over $X$, and $\a: E \to F$ a  morphism such that
    ${\a}_x: E_x \to F_x$ is an epimorphism for all points
    $x \in X$. Then there exists a morphism $\b: F \to E$,
    such that $\a \b = \Id_F$.
\end{th}

\dfi
Consider a strongly continuous action of $G$ on an arbitrary
     Banach space $\G$. A vector $s\in \G$ is called
    {\em periodic,}%
\label{ind:periodvec}
 if the orbit $Gs$ lies in a finite-dimensional
     subspace of space $\G$.

According to a lemma of Mostow (see \cite{Most}) periodic vectors form
a dense subset of $\G$.

\begin{th}\label{stabbund}
   Let $X$ be a compact $G$-space,
and $E \to X$ a $GGA$-bundle. Then there exists a
trivial $GGA$-bundle
   $\underline {M} =X \times M$ and $GGA$-bundle $E'$ such  that
   $\underline {M} \cong E \bigoplus E'$, where $M$ is a projective
$GGA$-module.
\end{th}

\Dok
By a lemma of Mostow \cite{Most} and Lemma~\cite{MF}
there exist sections $ \overline  s\,_1,\dots, \overline  s\,_n$
 periodic in $\G (E) $, as in complex Banach space, with
$\{\,\overline s_j(x)\}$ generating $E|_x$  for any $x$.
Taking their orbits,
  we consider finite-dimensional $G$-$\C$-module $W\ss\G (E) $
spanned by them: $W$ is the $\C$-linear span of
$\{G\,\overline s\,_j\}_{j=1,\dots,n}$.
Let $s_1,\dots, s_N$ be a basis for
 $W$. We have the trivial $GGA$-bundle
      $$
     \underline M = X \times M = X \times (A \o_{\C} W)
      $$
with the diagonal action of $G$. Define a morphism
$\te: \underline M \to E$ on the generators
by the formula $\te (x, a \o s_i) =s_i(x)\cdot a$.
This is an epimorphism of $A$-bundles. Since
      $$
  \te (g(x,a \o s_i))=\te(gx,ga \o g(s_i))=(g(s_i))(gx)\cdot ga=
      $$
      $$
   = gs_i(g_{-1}gx)\cdot ga =g(s_i(x))\cdot ga=g(s_i(x)\cdot a)=
g(\te(x,a\o s_i)),
      $$
   it follows that
   $\te$ is $G$-mapping. We introduce in $\underline M$
   the following fiberwise $A$-Hermitian product. If $(\,,\,)$ is a
   $G$-invariant inner product in $W$, then let
      $$
     \langle (x,a \o s_i),(x,b \o s_j)\rangle =a^*b(s_i,s_j).
      $$
   Then
\begin{eqnarray*}
\langle g(x,a \o s_i),g(x,b \o s_j)\rangle &=&g(a)^*g(b)(g(s_i),g(s_j))
=g(a^*b)(s_i,s_j)= \\
&=& g(\langle (x,a \o s_i),(x,b \o s_j)\rangle).
\end{eqnarray*}
Hence, our $A$-product is $G$-invariant (in the sense of Def.
\ref{def:invproizv}).
Such a product gives, in
particular, the structure of a Hilbert module in any fiber.
 Let $E'=\Ker \te$. With the help of \cite{MiUMN} it is
 easy to obtain that $\underline M \cong E ' \bigoplus E$.
\hfill$\Box$

 \begin{cor}
    A fiberwise $GGA$-product can be introduced on every
    $GGA$-bundle over a compact base.
 \end{cor}

   Let $E \to X$ be an $A$-bundle. We consider the bundle
$\Hom(E,E)$ with fiber $\Hom(E_x,E_x),\, x\in X$. An element
$p_x= {p_x}^2$ is called a {\em projection in the fiber}.
Denote by $\Proj(E) \ss \Hom(E, E)$ the set of projections.
Let $Q (E) \ss\Hom (E, E)$ consist of all $T$ such that
$z\cdot 1_{E_x} - T_x$ is an isomorphism for any
$x \in X,\, z \in \C$ if ${\rm Re}\:z\ne 1/2$.
It is clear, that $\Proj (E) \ss Q (E) $.
The usual method of Cauchy integrals (cf.~\cite[p. 184]{KarCli})
can be used to prove

 \begin{lem}
  There exists a retraction $\a: Q (E) \to \Proj (E) $.
  \hfill$\Box$
 \end{lem}

\vspace{-2ex}
\cst\label{conext} If $E'$ is a
  $GGA$-bundle over a closed $G$-invariant subspace $Y$ of a
  compact $G$-space $X$, then
  $E '\bop F' = Y \times M$ for some $F'$,
by Theorem~\ref{stabbund}. We define
  $p': Y \times M \to Y \times M$ by the formula
$p'(e',f')=(e',0)$, bearing in mind the
 identifications indicated above.
Then $p'$ is a
  $G$-projection with $\Im p'=E'$. We extend $p'$ as a cross-section
in $\G^G (\Hom (E, F))$ over $Y$ to a $G$-section $p$ over $X$, where
  $F=X \times M$. By Lemma on generators from \cite{MF},
there exists a neighborhood
  $U$ of $Y$ such that
  $p\vert_{\overline U} \in Q(F\vert_{\overline U})$.
  Taking $GU$, we may assume that $U$ is
  $G$-stable. It is easy to see that $E=\Im\a(p\vert_{\overline U})$
 is a $GGA$-bundle and $E\vert_Y = E'$.

 \dfi
   Following the general scheme presented above, we define the
   $K$-groups for a compact $G$-space $X$ by setting
      $$
     \begin{array}{rl}
     \bK^{p,q}_G(X,A)&=K^{p,q}(\bVect _G(X,A)),\\[5pt]
     \bK^{n}_G(X,A)&=K^{n}(\bVect _G(X,A)),\\[5pt]
     \bK_G(X,A)&=K^{0,0}(\bVect _G(X,A)).
     \end{array}
      $$

Let $\r^{X,Y}: \bVect _G(X,A) \to \bVect _G(Y,A)$
be the restriction functor, where $Y$ is a closed $G$-subspace of
  $X$. The following assertion is a consequence of
  Theorem~\ref{stabbund} and Lemma~\ref{prodoljmor}.

 \begin{lem}
$\r^{X,Y}$ is a full quasi surjective
Banach functor in the sense of {\rm~\ref {kar4}} .\hfill$\Box$
 \end{lem}

Therefore, setting $\bK^n_G(X,Y,A)  :=  \bK^n(\ro^{X,Y})$ in the sense
of Definition \ref{kar9a}, we get (by \cite[II.3.22]{KaKn} and
\cite[2.3.1]{KarCli}) the exact sequence of a pair
   $$
   \dots \to \bK^{n-1}_G(X,Y,A)\to \bK^{n-1}_G(X,A)\to
\bK^{n-1}_G(Y,A)\to
   $$
   \begin{equation}\label{exact}
\to
   \bK^n_G(X,Y,A) \to \bK^n_G(X,A) \to \bK^n_G(Y,A) \to \dots
    \end{equation}

Consider the trivial $G$-pair $(B^n, S^{n-1})$, where $B^n$
is the $n$-dimensional closed ball, and $S^{n-1}$ its
boundary.  Then, according to a construction in~\cite{KaKn,KarCli}
for an
arbitrary Banach category, the categories
$\bVect_G(X,A)(B^n)$, $\bVect_G(X,A)(S^{n-1})$
are defined along with the group
      $
     \bK(B^n,S^{n-1};\Vect_G(X,A))=K(\vfi),
      $
where
      $$
     \vfi : \bVect_G(X,A)(B^n) \to \bVect_G(X,A)(S^{n-1})
     $$
is the restriction functor.

We now give the definitions and show, that $\f$ coincides with
      $$
     \ro^{B^n \times X, S^{n-1}\times X}:\bVect_G(B^n \times X,A)
     \to \bVect_G(S^{n-1} \times X,A).
      $$

Since both functors are induced by restrictions, it is suffices
to show, that for a $G$-trivial compact space $Z$ the categories
$\bVect_G (Z \times X, A)$ and $\bVect_G (X, A) (Z)$ are naturally
isomorphic to each other. By definition, $\bVect_G (X, A) (Z)$ is
associated in the sense of \ref{kar1} with the category
$\bVect_G (X, A)_T (Z)$, defined as follows. The objects of
$\bVect_G (X, A)_T(Z)$ coincide with those of $\bVect_G(X,A)$, and
the morphisms are continuous mappings
$$
Z \to \Mor_{\bVect_G(X,A)}(E,F).
$$
Thus, $\bVect_G (X,A)_T (Z) $ is identified with a full subcategory in
$\bVect_G(X\times Z, A)$. The category
$$
\bVect^T_G (X \times Z,A)
$$
of trivial $GGA$-bundles over $X \times Z$, with which
$\bVect_G (X \times Z, A)$ is associated according to Theorem
\ref{stabbund} is a full subcategory of $\bVect_G (X, A)_T (Z)$.
Passing to the associated categories, we get the required
result. Comparison of it with constructions in \cite{KaKn,KarCli},
yields two important corollaries.

\begin{th}\label{bott} {\rm (Bott-Clifford periodicity) \/}
We have a natural isomorphism
      $$
     \bK^n_G(X,Y,A) \cong \bK^{n-2}_G(X,Y,A).
      $$
 \end{th}
 \begin{th} {\rm (Bott periodicity) \/}
We have a natural isomorphism
      $$
     \bK^1_G(X,A) \cong \bK_G(X\times B^1,X\times S^0,A).
      $$
  \end{th}
We get the next result from Construction~\ref{conext}.
   \begin{lem}\label{limotkr}
Let $Y$ be a closed $G$-subspace of a compact space $X$.
Then
    $$
   \bK_G(Y,A) = \lim\limits_{\ola{Y \ss U}} \bK_G(U,A).
   $$
   \end{lem}

   \dfi
Let $X$ be a locally compact paracompact (Hausdorff)
$G$-space.  Let
      $$
     \bK_{G,c}(X,A) := \Ker \{\bK_G(\dot X,A) \to \bK_G(pt,A)\},
      $$
where $\dot X = X \cup pt $ is the one-point
compactification.

It is easy to get another description of groups $K_{G, c} $
from a consideration of the sequence (\ref{exact}) for the pair
$(\dot X, pt)$ :
 \begin{lem}
      $$
     \bK_{G,c}(X,A) \cong \bK_G(\dot X, pt, A).
      $$
 \end{lem}

\dfi
Let
      $$
\bK^{-1}_{G,c}(X,A)=\Ker\{\bK^{-1}_G(\dot X, A) \to \bK^{-1}_G(pt,A)\},
      $$
      $$
\bK^n_{G,c}(X,A)=\bK_{G,c}((X \setminus Y) \times \R ^n,A).
      $$

It is necessary to verify compatibility:
      $$
     \bK^{-1}_{G,c}(Y,A) \cong \bK_{G,c}(Y \times \R, A).
      $$
This is done in the same way as for Theorem \cite[II.4.8]{KaKn}. If $X$
and $Y$ are compact, then
$\bK^n_{G,c}(X,Y,A)\cong \bK^n_G(X,Y,A)$;
therefore, we omit the sign $c$.

We finish this Section with another description of our K-groups.

\dfi Let $X$ be a locally compact paracompact
$G$-space. A {\em complex of $GGA$-bundles\/}
over $X$ is defined to be a sequence
 $$
(E,d)=\Bigl(\dots\ola{d_i}E^i\ola{d_{i+1}}E^{i+1}\ola{d_{i+2}}
\dots\Bigr),\qquad i\in\Z,
 $$
where $i\in\Z$, the $E^i$ are $GGA$-bundles over $X,$
and $d_i$ are morphisms with $d_ {i + 1} d_i=0$ for every $i$;
also, $E^i=0$ for all but perhaps finitely many $i$.
A {\em morphism of complexes\/} $f:(E,d)\to (F,h)$
is defined to be a sequence of morphisms $f_i: E^i\to F^i$ connected
by the condition $f_{i+1} d_{i+1} =h_{i+1}f_i.$ Isomorphism in
this category will be denoted by $(E, d)\cong (F, h).$
A point $x\in X$ is called a {\em point of acyclicity $(E,d)$\/}
if the restriction of $(E,d)$ to $x,$ i.~e., the sequence of $A$-modules
$$
(E,d)_x=\Bigl(\dots\ola{(d_i)_x}E^i_x\ola{(d_{i+1})_x}E^{i+1}_x
\ola{(d_{i+2})_x}\dots\Bigr),
 $$
is exact. The {\em support $\supp (E, d)$\/}
is the complement in $X$ of the set of points of acyclicity.

 \begin{prop}
$\supp (E, d) $ is a closed $G$-subset of $X.$
 \end{prop}

\Dok
The $G$-invariance is obvious. The proof of the fact that it is closed
can be found, e.~g., in \cite[1.3.34]{SolTro}.
\hfill$\Box$

\rem
The remaining assertions of this section can be proved by the classical
scheme (cf.~\cite{AK,Fried}) with the specific nature of C*-algebras
taken into account, as it was demonstrated in two previous lemmas;
therefore, the proofs are omitted.

\dfi
Let $Y$ be a closed $G$-subspace of $X.$ Denote
by $\bL_G (X, Y; A) $ the semigroup (with respect to the direct sum)
of classes of
isomorphic $GGA$-complexes $(E,d)$ on $X$ such that $\supp (E,d)$ is
a compact subset of $X\setminus Y.$ Two elements in $\bL_G (X, Y; A) $
are said to be {\em homotopic\/} if for representatives
$(E_0,d_0)$ and $(E_1,d_1)$ of them there is a complex
$(E,d)$ representing an element of
$\bL_G (X\times I, Y\times I; A)$ such that
$(E_0, d_0)= (E, d)|_{X\times \{0\}}$
and $(E_1,d_1) = (E, d)|_{X\times \{1\}}$;
in this case we write
$(E_0, d_0) \simeq (E_1, d_1).$ We introduce the following
equivalence relation: $(E_0,d_0)\sim (E_1, d_1)$ iff there are
acyclic $(F_0, f_0)$ and $(F_1, f_1)$ such that
 $$
(E_0,d_0)\oplus (F_0,f_0) \simeq (E_1,d_1)\oplus (F_1,f_1).
 $$
By
acyclicity we understand the condition $\supp (E, d)=\emptyset.$
Let
$
\bM_G (X, Y; A) =\bL_G (X, Y; A) /\sim.
$
Denote by $\bL_G^n (X, Y; A) $ and $\bM_G^n (X, Y; A)$ the corresponding
semigroups constructed from the complexes of length $n.$
We have the natural injective semigroup homomorphisms
(addition of the zero term):
 $$
\bL_G^n(X,Y;A)\to \bL_G^{n+1}(X,Y;A),\quad \bL_G(X,Y;A):=\lim_\to
\bL_G^n(X,Y;A).
 $$
The equivalence relation $\sim$ commutes with imbeddings;
therefore, the indicated morphisms induce morphisms
 $
\bM_G^n(X,Y;A)\to \bM_G^{n+1}(X,Y;A).
 $
 \begin{lem}
Suppose that $E\to X$ and $F\to X$ are
$GGA$-bundles, $\a: E|_Y\to F|_Y$ and $\b: E\to F$ are
monomorphisms, and in the class of monomorphisms
$E|_Y\to F|_Y$ there exists a $G$-homotopy joining $\a$ and $\b|_Y.$
Then there exists a monomorphism $\widetilde\a: E\to F$ such
that $\widetilde\a|_Y=\a.$ \hfill$\Box$
\end{lem}

\begin{lem}
$\bM_G^n(X,Y;A)\to \bM_G^{n+1}(X,Y;A)$ is an isomorphism.
\hfill$\Box$
 \end{lem}

\rem\label{rem:ka25}
Suppose that $X$ is compact and $Y=\emptyset.$ Then we have a natural
isomorphism $\x_1:\bM_G^1(X,\emptyset;A)\to \bK_G(X;A),$ assigning
to a class of the complex $(0\to E^0\ola {\a} E^1\to 0)$ the element
$ [E^1] - [E^0].$

\begin{lem}
There exists and unique natural
equivalence of functors
 $$
\x_1:\bM_G^1(X,Y;A)\to \bK_G(X,Y;A),
 $$
on the category of compact $G$-pairs, with the equivalence
of the form indicated in {\rm \ref{rem:ka25}}
for $(X,\emptyset)$.
 \end{lem}

\dfi
We return to the case of (locally
compact, paracompact) noncompact, in general, $G$-pairs
$(X,Y)$. Two complexes $(E, d)$ and $(F,b)$ in $\bL_G (X, Y; A)$
are said to be {\em compactly isomorphic\/}
\label{ind:gabundlecomi}
if there exists a compact $G$-subset $C$ with
 $$
\supp (E,d)\cup\supp (F,b)\subset {\rm int}\,C\subset C\subset
X\setminus Y,
 $$
and $G$-isomorphisms $\psi^i: E^i\to F^i$ over $X,$
where the following diagram commutes over
$X\setminus {\rm int} \, C$:
 $$
\begin{array}{ccccccccccc}
0 & \ola{} & E^0 & \ola{d^1} & E^1 & \ola{d^2} & \dots & \ola{d^n} &
E^n &\ola{}& 0 \\
 & & \quad\downarrow\,\psi^0 & & \quad\downarrow\,\psi^1 & & & &
\quad\downarrow\,\psi^0 & \\
0 & \ola{} & F^0 & \ola{b^1} & F^1 & \ola{b^2} & \dots & \ola{b^n} &
F^n &\ola{} & 0.
\end{array}
 $$

Complex of the form
$0\to\dots\to 0\to E\ola {\Id} E\to 0\to \dots \to 0$
Is said to be {\em elementary.} Two complexes $(E, d)$ and $(F, b)$
are said to be {\em compact equivalent,\/}
\label{ind:gabundlcomcompequi}
if there exist elementary $Q_1,\dots, Q_r$ and $P_1,\dots, P_s$ such
that
$ (E, d) \oplus Q_1 \oplus\dots\oplus Q_r$ is compact isomorphic to
$ (F, b) \oplus P_1 \oplus\dots\oplus P_s.$ This equivalence is
denoted by $\approx.$

 \begin{lem}
$(E,d)\approx (F,b)$ implies $(E,d)\sim (F,b).$ \hfill$\Box$
 \end{lem}

 \begin{lem}
Suppose that
$$
(E,d)=\Bigl(0\to E^0\ola {d} E^1\to 0\Bigr),\qquad
(F, b) =\Bigl (0\to F^0\ola {b} F^1\to 0\Bigr)
$$
and there exist compact $G$-sets $C_1$ and $C_2$ such
that
 $$
\supp (E,d) \cup \supp (F,b) \subset C_1 \subset {\rm int}\,C_2
\subset C_2 \subset X\setminus Y,
 $$
and a $GGA$-bundle $L$ over $C_2$ with an isomorphism
 $$
\te : (E^0\oplus F^1)|_{C_2}\oplus L \to (E^1\oplus F^0)|_{C_2}\oplus
L,
 $$
for which $\te|_{C_2\setminus C_1}=d\oplus b^{-1}\oplus 1.$
Then
$ (E, d) \approx (F, b).$\hfill$\Box$
\end{lem}

 \begin{lem}
$\bM_G^1(X\setminus Y;A)\cong \bM_G^1(X,Y;A).$ \hfill$\Box$
 \end{lem}

 \begin{cor}
 $$
\begin{array}{llcl}
1)\quad & \bM_G^1(\dot X, pt; A) & \cong & \bM_G^1(X;A)    \\
        & \downarrow\cong     & & \\
        & \bK_G(\dot X, pt; A) & \cong & \bK_G(X;A),    \\[15pt]
2)\quad & \bM_G^1(X,Y; A) & \cong & \bM_G^1(X\setminus Y;A)    \\
        & & &\downarrow\cong      \\
        & \bK_G(X,Y; A) & \cong & \bK_G(X\setminus Y;A)    \\
\end{array}
$$
\hfill$\Box$
 \end{cor}

 \begin{th}\label{KiM}
There is a natural isomorphism
 $$
\bM_G(X,Y;A)\cong \bM_G^1(X,Y;A)\cong \bK_G(X,Y;A). \qquad\Box
 $$
 \end{th}

\section{С*-Fredholm operators}

Let us remind the definition of C*-Fredholm operator \cite{MF,TroIAN}.

\dfi\label{fredbezg}
A bounded $A$-operator $F: H_A \to H_A$,
is called {\em Fredholm}, if

1) The operator $F$ supposes adjoint, and

2) There exist
decompositions of the domain $H_A=M_1\bop N_1$ and
   the values space
   $H_A=M_2\bop N_2$ of $F$ (where $M_1, M_2, N_1, N_2$ are closed
   $A$-submodules, $N_1$ and $N_2$ have a finite number of
generators) such that
   the operator $F$ has with respect to these decompositions
the matrix
   $F=\Mat {F_1}00{F_2}$, and $F_1: M_1\to M_2$
   is an isomorphism.

In the equivariant case we change the definition as follows.
Let $l_2(P)$ be equipped with an invariant $A$-inner product
(it is possible to do this by Remark~\ref{usrP}).
We can apply  the Stabilization Theorem of Kasparov \cite{KaspJO}
to $l_2 (P)$:
\begin{equation}\label{fred*}
l_2(P)\oplus\H_A \cong \H_A,
\end{equation}
where  $\H_A=\sum_{i=1}^\infty(A\o_\C V_i)$, and $\{V_i\}$
is the countable collection of
finite-dimensional spaces, in which
    all (up to isomorphism) irreducible unitary
     representations of $G$ are  realized, and each representation
is repeated an infinite number of times.
Isomorphism (\ref{fred*}) is a
     $GGA$-isomorphism of Hilbert modules and the sum on the left
in (\ref{fred*}) is orthogonal. We introduce the following notation
 $$
R_m=\sum_{i=m+1}^\infty (A\o_\C V_i),\qquad
R_m^\bot=\sum_{i=1}^m (A\o_\C V_i).
      $$

     For any bounded $G$-$A$-operator $F: l_2 (P_1) \to l_2 (P_2) $
let $S (F):\H_A\to\H_A$
     ($S$ due to the word ``stabilization'')
     denote the operator
 $$
\Mat 100F : \H_A\cong \H_A \oplus l_2(P_1) \to \H_A \oplus l_2(P_2)
\cong \H_A.
      $$
Of course, everything is determined up to a
$GGA$-isomorphism.

\begin{th}\label{fred1}{\rm cf.~\cite{MF,TroIAN} \/}
  Let $\H_A\cong \M\oplus \N$,
  where $\M$ and $\N$ are closed $GGA$-modules, and $\N$ has
  a finitely many generators $a_1,\dots, a_s$. Then $\N$ is a
  projective $GGA$-module of finite type.
\end{th}

\Dok Just as in~\cite{MF,TroIAN}.\hfill$\Box$

\dfi\label{fredopr}
A bounded $GGA$-operator
     $$
     F:l_2(\P_1)\to l_2(\P_2),
     $$
is called {\em Fredholm operator\/} ($GGA$-Fredholm), if

1) $F$ admits an adjoint;

2) for $S(F)$ there exist an inverse image decomposition
   $\H_A=\M_1\bop\N_1$ and an image decomposition
   $\H_A=\M_2\bop\N_2$, where $\M_1,\,\M_2,\,\N_1,\,\N_2$ are closed
   $GGA$-modules, $\N_1,\,\N_2$ have finitely many generators,
   and the operator $S(F)$ has the matrix form
   $S(F)=\Mat {F_1}00{F_2}$
   in these decompositions, with $F_1:\,M_1\to \,M_2$ being
   an isomorphism of $GGA$-modules. By Theorem~\ref{fred1},
   $N_1$ and $N_2$ are projective $GGA$-modules; we can form
   the index element
    $$
\index F=[\N_1]-[\N_2]\in \bK^G(A).
    $$

 \begin{th}\label{horfred}
{\rm (see, e.~g.~\cite{SolTro})}
In the decomposition in the definition of $A$-Fredholm
     operator
      {\em ({\em see {\em \ref{fredbezg} \/} \/}) \/}
we can always assume $M_0$ and $M_1$ admitting an
     orthogonal complement. More precisely, there exists a
     decomposition for  $F$
       $$
     \Mat{F_3}00{F_4} :H_A=V_0\bop W_0 \to
      V_1\bop W_1 = H_A,
      $$
     such that $V^\bot_0 \hop V_0 =H _A,\,V^\bot_1\hop V_1=H_A$,
     or (what is just the same) such that the
     projections
      $$
     p_0:V_0 \bop W_0 \to V_1, \quad p_1:V_1\bop W_1 \to V_1
      $$
     admit conjugates.
 \end{th}

 \begin{lem}
Let $\M=\M_0\oplus \M_1$ be a decomposition into an orthogonal sum and
$\M$ and $\M_0$ are $GGA$-modules with an invariant inner product.
Then $\M_1$ is invariant.
 \end{lem}

\Dok Let $x\in\M_0$, $y\in\M_1$ and $g\in G$ be arbitrary. Then
$$
 \<x,gy\>= g(\<g^{-1}x,y\>)=g(0)=0.\qquad \Box
$$

 \begin{cor}
The previous Theorem is valid in the $GGA$-case as well.
 \end{cor}

\Dok
While proving the mentioned theorem we took
$V_0=(N_0)^\bot$ and $V_1=F(V_0)$. Hence, by Lemma all modules are
invariant.\hfill $\Box$

The remaining statements of the present section can be proved
similarly to the $GA$-case (see~\cite{TroIAN,SolTro}).

 \begin{th}\label{corrind}
The index is well-defined.
 \end{th}

  \begin{lem}\label{fred6}
Let $F:l_2 (A)\to l_2(A)$ be a
$GGA$-Fredholm operator,
then there exists an $\e > 0$ such that any bounded
 $GGA$-operator $D$,
 restricted to satisfy $\| F-D\| < \e$ and
 admitting an adjoint,
is a $GGA$-Fredholm operator and $\index D=\index F$.
  \end{lem}

  \begin{lem}\label{compozfred}
 Let $F$ and $D$ be $GGA$-Fredholm operators,
 $$
F:l_2(\P_1)\to l_2(\P_2), \qquad D:l_2(\P_2)\to l_2(\P_3).
  $$
 Then $DF: l_2 (\P_1) \to l_2 (\P_3) $ is a $GGA$-Fredholm
 operator and $\index DF=\index D + \index F$.
  \end{lem}

 \begin{lem}\label{fred9}
Let $K: l_2 (\P) \to l_2 (\P) $ be a compact $GGA$-operator.
Then $1+K$ is a $GGA$-Fredholm operator and
$\index (1 + K) =0$.
 \end{lem}

\begin{lem}\label{fred10}
Consider a $GGA$-Fredholm operator $F:l_2(\P_1) \to l_2(\P_2)$.
Let an operator $K\in\K (l_2 (\P_1), l_2 (\P_2))$ be
$G$-equivariant. Then the operator $F + K$ is a
$GGA$-Fredholm operator and $\index (F + K) =\index F$.
 \end{lem}

  \begin{lem}\label{lem:fredh62}
Let $F: l_2 (\P_1) \to l_2 (\P_2) $ be a bounded $GGA$-operator
admitting an adjoint, $D=S (F) + K\in\End^* \H_A$ and
$K\in\K (\H_A)$ be $GGA$-operators. Let $D$ have a decomposition of
$\H_A$ from the definition of $G$-$A$-Fredholm operator.
Then $F$ is a $GGA$-Fredholm operator.
 \end{lem}

 \begin{th}\label{kritfred}
Let
 $$
F:l_2(\P_1)\to l_2(\P_2),\quad D:l_2(\P_2)\to l_2(\P_1),\quad
D':l_2(\P_2)\to l_2(\P_1)
 $$
be bounded $GGA$-operators admitting an adjoint and
 $$
S(FD)=1_{\H_A}+K_1,\quad S(D'F)=1_{\H_A}+K_2,\quad K_1,\,
 K_2\in\K(\H_A).
 $$
 Then $F$ is a $GGA$-Fredholm operator.
 \end{th}

 \begin{lem}
Let $D$, $D'$, and $F$ be bounded $GGA$-operators admitting
adjoint. Let $FD$ and $D'F$ be $GGA$-Fredholm
operators.Then $F$ is a $GGA$-Fredholm operator.
 \end{lem}

\section{The Thom isomorphism}\label{subs:Thom}

In this section we discuss the theorem on
the Thom isomorphism in $\bK_G(\:,A)$-theory. As in other cases
(see e.g. \cite{FaSk}) it plays an important role.

Let $X$ be a $G$-space, $p: F\to X$ a complex bundle over $X$, and
$s: X\to F$ an invariant section. We denote  by $\La^i(F)$ the
complex $G$-bundle of $i$-vectors over $X$. Let us define
the complex $\La (F,s)$ over $X$ of lengths $n=\dim F$:
 $$
\La(F,s):=(0\to
\La^0(F)\ola{\a^0}\La^1(F)\ola{\a^1}\dots\ola{\a^{n-1}}\La^n(F)\to 0),
 $$
where $\a^k(v_x)=s(x)\w v_x$ for $v_x\in\La^k(F)_x$.

\begin{lem}
$($ see {\rm \cite{Fried})}

\noindent
{\rm\bf 1.} $(\La (F,s),\a)$ is really a complex.

\smallskip
\noindent
{\rm\bf 2.} $\supp (\La(F,s))=\{x\in X|s(x)=0\}$ and if
this set is compact, then the element
$[\La(F,s)]\in K_G(X)$ is defined.\hfill$\Box$
\end{lem}

Let $\pi: p^* F\to F$ be the bundle with total space
 $$
p^*F=\{(f_1,f_2)\in F\times F\: |\: p(f_1)=p(f_2)\},
\qquad \pi (f_1,f_2)=f_1.
 $$
The vector bundle $(p^*F,\pi, F)$ has the canonical section
 $$
s_F:F\to p^*F,\qquad s_F(f)=(f,f).
 $$
The support of $s_F$ is equal to $X$. Hence, if $X$ is a compact set,
the element
$
 [\La(p^*F,s_F)]=\la_F\in K_G(X).
$
is defined. Let us define by $a\cdot b$ the element
$p^*(a)\o b \in \bK_G (F;A).$

\dfi
If the base of a vector bundle is compact, then
the mapping
$$
\f:\bK_G(X;A)\to \bK_G(F;A),\qquad \f(a)=a\cdot \la_F.
$$
is called the
{\em Thom homomorphism.\/}
\label{ind:Thomhomomorphism}

The following statement is obvious.
\begin{prop}\label{prop:thom3}
The Thom homomorphism
$\f$ is a morphism of $R(G)$-mo\-du\-les. \hfill$\Box$
\end{prop}

Let $i: X\hookrightarrow F$ be the enclosure of the zero section
of $F$. It induces the homomorphisms
 $$
i^*:\bK_G(F;A)\to \bK_G(X;A),\qquad i^*\f:K_G(X;A)\to K_G(X;A).
 $$

\begin{prop}
{\rm\bf 1.} If $X$ is compact, then the following sequence is defined
 $$
0\lra\La^0F\ola 0 \La^1F\ola 0 \La^2F\lra\dots\ola 0 \La^nF\lra 0.
 $$
It defines an element of $K_G(X)$. For any element
$a\in \bK_G(X;A)$
$$
i^*\f(a)=a\cdot[0\lra \La^0F\ola 0 \dots \ola 0 \La^n F\lra 0].
$$

\smallskip
\noindent
{\rm\bf 2.} For $a\in \bK_G(X;A)\cong \bM_G(X;A)$ $($see {\rm\ref{KiM})}
$$
i^*\f(a)=a\cdot\sum_{i=0}^n(-1)^i\La^iF.
$$
 \end{prop}

\Dok
By the construction of the natural isomorphism between $\bK_G$ and
$\bM_G$, these assertions are equivalent to each other.
Let us prove the item 1. Let $a\in \bK_G(X;A)$
be represented by the complex
 $$
a=[(\E,\a)]=[0\lra E^k\ola {\a^k}\dots \ola {\a^1}E^0 \lra 0].
 $$
Then
 $$
\f(a)=a\cdot\la_F=[(p^*\E,p^*\a)\o\La(p^*F,s_F)],
 $$
hence
 $$
i^*\f(a)=[(i^*p^*\E,i^*p^*\a)]\o [i^*\La(p^*F,s_F)].
 $$
Since $pi=1$, we have
$(i^*p^*\E,i^*p^*\a)=(\E,\a)=a$, while the restriction of complex
$\La(p^*F,s_F)$ on $X$ is equal to
 $$
[0\lra \La^0F\ola 0 \dots \ola 0 \La^n F\lra 0].
\qquad\Box
 $$

Let us pass to the case of a locally compact space $X$.
The complex $\La (p^* F, s_F)$
has no compact support now and does not determine the
element $\la_F\in K_G (F)$. However, if $a= [(\E,\a)] \in \bK_G (X; A)$,
then
\begin{eqnarray*}
\supp \{(p^*\E,p^*\a)\o &&\hspace{-1.7em}
\La(p^*F,s_F)\}\ss\\
&&\ss \supp (p^*\E,p^*\a)\cap\supp\La(p^*F,s_F)\ss\\
&&\ss\supp (p^*\E,p^*\a)\cap X=\supp (\E,\a).
\end{eqnarray*}
Thus, the complex $(p^*\E,p^*\a)\o \La(p^*F,s_F)$ has the compact
support.
We obtain a homomorphism of $R(G)$-modules
 $$
\f:\bK_G(X;A)\to \bK_G(F;A),\qquad \f(a)=[(p^*\E,p^*\a)\o \La(p^*F,s_F)].
 $$
As well as in the compact case,
 $$
i^*\f(a)=a\,[0\to\La^0F\ola{0}\La^1F\ola{0}\dots\ola{0}\La^nF\to 0].
 $$
Passing by the Bott periodicity to $K^1$ (see Theorem~\ref{bott}),
we define the {\em Thom homomorphism\/}
\label{ind:Thomhomomorphism2}
in the general case:
 $$
\f=\f^F_A:\bK^*_G(X;A)\to \bK^*_G(F;A).
 $$

Let now $E$ and $F$ be two complex $G$-bundles over $X.$
Then a product is defined by the formula
 $$
\mu:K_G(E)\o K_G(F)\to K_G(E\times F).
 $$
Since the enclosure $E\oplus F \to E\times F$ induces a homomorphism
 $$
K_G(E\times F)\to K_G(E\oplus F),
 $$
we obtain the multiplication
 $
\cdot : K_G(E)\o K_G(F) \to K_G(E\oplus F).
 $

 \begin{th}\label{svolamb}
{\rm\bf 1.} {\rm \cite{Fried}} If $X$ is compact, then
     $$
     \la_E\cdot \la_F = \la_{E\oplus F}.
      $$

\smallskip
\noindent
{\rm\bf 2.} In the general case
 $
\f^E_A\cdot\f^F_\C=\f^E_\C \cdot \f^F_A= \f^{E\oplus F}_A.
 $
The equalities hold
in the following sense. Let, for example, $ab\in \bK_G(X;A),$ where
$a\in  K_G(X),$ $b\in \bK_G(X;A)$. Then
      $
     \f^E_\C (a)\cdot \f^F_A(b)= \f^{E\oplus F}_A(ab).
      $

\smallskip
\noindent
{\rm\bf 3.} {\rm \cite{Fried}}
Let $F_1$ and $F_2$ be two complex bundles
over $X,$ and $s_1$, $s_2$ their section. Then
 $$
\La(F_1\oplus F_2, s_1\oplus s_2)=\La(F_1,s_1)\o \La(F_2,s_2).
 $$
 \end{th}

\Dok
The item 2 immediately follows from the item 1 and
the definition of the Thom homomorphism. The items 1 and 3 are proved
in~\cite{Fried}. \hfill$\Box$

 \begin{prop}\label{prop:thom6}
Consider $E\oplus F$ as a bundle over $X$ as well as a bundle
over $E.$ Then the diagram
 $$
\begin{array}{ccc}
\bK^*_G(X;A) & \ola{\f^E} & \bK^*_G(E;A) \\
\quad\downarrow{\scriptstyle \f_1^{E\oplus F}}&
& \quad\downarrow{\scriptstyle \f_2^{E\oplus F}}\\
\bK^*_G(E\oplus F;A) & = & \bK^*_G(E\oplus F;A)
\end{array}
 $$
is commutative.
 \end{prop}

\Dok
Consider the projections
 $$
p:E\to X,\quad q:F\to X, \quad r:E\oplus F \to X,\quad t: E\oplus F
\to E.
 $$
Let $x\in K_G(X;A).$ Then
 $
\f^E(x)= p^*(x)\,\La(p^*E,s_E),
 $
 $$
\f^{E\oplus F}_2 \f^E(x)=t^*\Bigl(p^*(x)\La(p^*E,s_E)\Bigr)\,
                   \La\Bigl(t^*(E\oplus F),s_{E\oplus F}\Bigr),
 $$
 $$
t^*p^*(x)=r^*(x),\qquad t^*\La(p^*E,s_E)=\La(r^*E,t^*s_E)
 $$
and
  $
\La(t^*(E\oplus F),s_{E\oplus F})= \La (r^*F,t^*s_F).
 $
Since $t^*s_E +t^*s_F=s_{E\oplus F},$ by the item 3 of
Theorem~\ref{svolamb} (we denote the elements of K-groups and
their representatives by the same symbols),
 $$
\f^{E\oplus F}_2 \f^E(x)=r^*(x)\La(r^*E,t^*s_E)\La(r^*F,t^*s_F)=
 $$
 $$
=r^*(x) \La(r^*(E\oplus F),s_{E\oplus F})=\f^{E\oplus F}_1 (x).
\qquad\Box
$$

Let starting from this moment $X$ be {\em separable.\/}
Then $C_0(X)\skr G$ is of the class of algebras for which the results
of \cite{SchoII} are valid.

 \begin{th}\label{isotoma}~{\rm\cite{TroStekl}}
Let $X$ be separable and metrizable, and let action of
$G$ on $A$ be trivial. Then $\varphi_A$ is an isomorphism.
 \end{th}

\Dok
Let us form the geometrical resolution \cite{SchoII}
 $$
0\to A\ot\K\ot C_0(\R)\ola{\io} C\ola{\nu}0,
 $$
where $C\ne\C$ and $F$ are $C^*$-algebras, and $\K$ is the algebra of
compact operators in Hilbert space.
Let us prove the following Lemma.

 \begin{lem}\label{lem:kommut}
The diagram
 $$
\diagram
\bK^*_G(X;A\ot\K\ot C_0(\R))\rto^{\qquad \io_*}
\dto^{\f_{A\ot\K\ot C_0(\R)}}
&\bK^*_G(X;C)\dto^{\f_C}\\
\bK^*_G(V;A\ot\K\ot C_0(\R))\rto^{\qquad \io_*}&\bK^*_G(V;C)
\enddiagram
 $$
is commutative. Here $V\to X$ is a $\C$-vector bundle,
 $$
\io: A\ot\K\ot C_0(\R)=\Ker\nu\hookrightarrow C,
 $$
$\io_*$ is an extension of scalars.
 \end{lem}

\Dok
Let $a\in \bK_G^0(X;A\ot\K\ot C_0(\R))$ be defined by the following
complex:
 $$
(\E,\a)=(0\to E^1\ola{\a} E^0\to 0).
 $$
Then (by $\wo$ we define the tensor product over the cartesian product
of bases)
 $$
\f_{A\ot\K\ot C_0(\R)}(a)=[(p^*\E,\p^*\a)\ot \La(p^*F,s_F)]=
 $$
 $$
\diagram
=\Bigl[ 0\rto& p^*E^1\wo p^*\La^0 F \qquad
\rto^{(p^*\a\wo 1)\oplus (1\wo s)\qquad\qquad  }
&\qquad
(p^*E^1\wo p^*\La^0 F)\oplus (p^*E^1\wo p^*\La^1 F)\rto&\dots
\Bigr] \:|_{\rm diag},\\
\enddiagram
 $$
and
 $$
\io_*\f_{A\ot\K\ot C_0(\R)}(a)=
\Bigl[0\to (p^*E^1\ot_{A\ot\K\ot C_0(\R)}C)\wo p^*\La^0 F
\ola{(p^*\a\ot 1_C\wo 1)\oplus ((1\ot 1_C)\wo s)\qquad }
 $$
 $$
\to ((p^* E^0\ot_{A\ot\K\ot C_0(\R)}C)\wo p^*\La^0 F)\oplus
((p^* E^1\ot_{A\ot\K\ot C_0(\R)}C)\wo p^*\La^1 F)\to\dots
\Bigr] \:|_{\rm diag}=
 $$
 $$
=\f_C\io_*(a).\qquad \qquad \qquad \qquad \qquad \qquad \qquad
 $$
By definition $\f$ commutes with the Bott isomorphism,
which is natural \cite{KarCli}, in particular, it commutes with
$\iota_*$. Hence, Lemma is true for $K^1$ too.
\hfill$\Box$

If $G$ is a compact group, $X$ is a paracompact Hausdorff
$G$-space, then $X/G$ is a paracompact
$G$-space (see, e.~g. \cite[p. 7]{Fried}). The stabilizer
$G_x$ is closed \cite{Bre}. Hence, it is compact and
has the type I. Therefore, by the following theorem,
C*-crossed product
$C_0(X)\skr G$ has the type I in our case.

 \begin{th}
{\rm~\cite{Goot}}
$C_0(X)\skr G$ has the type {\rm I} iff
$X/G$ is a $\tau_0$-space and all isotropy groups
$G_x$ has the type {\rm I}. \hfill$\Box$
 \end{th}

Let $\g:K^*_G(X;B)\ola\cong K_*((C_0(X)\ot B)\skr G)$
be the natural isomorphism \cite{Green,Julg}.
Let $B:=A\ot\K\ot C_0(\R)$ and $C^X_G:=C_0(X)\skr G$ (and similarly
for $V$). Consider the following diagram
 $$
\diagram
0\to\Tor(K_*(C^X_G),K_*(A))\rto
\dto_\cong^{\g\f_\C\g^{-1}*1\qquad \fbox{1}}
& K_*(C^X_G)\ot K_*(C)\rto\dto_\cong^{\g\f_\C\g^{-1}\ot 1\qquad
\fbox{2}}
&K_*(C^X_G\ot C)\dto^{\g\f_C\g^{-1}}\\
0\to\Tor(K_*(C^V_G),K_*(A))\rto
& K_*(C^V_G)\ot K_*(C)\rto
&K_*(C^V_G\ot C),
\enddiagram
 $$
where the exact rows represent a part of (4.5) from \cite{SchoII}.

Consider also the diagram
 $$
\diagram
&\bigoplus\limits_{q\in\Z_2}K_q(C^X_G\ot B)
\rto^{\iota_*}
& \bigoplus\limits_{q\in\Z_2}K_q(C^X_G\ot C)
&\\
0\to K_*(C^X_G)\ot K_*(A)\rto^{\qquad \a_X}
&K_*(C^X_G\ot A)\rto^{\b_X\qquad }
\uto_{\cong\qquad \qquad \fbox{8}}^{\rm Bott}
&\Tor (K_*(C^X_G), K_*(A))\to 0 \uto^{\bigcup\quad}
\\
&K_*((C_0(X)\ot A)\skr G)\uto_\cong^{\s_X}&
\ddto^{\g\f_\C\g^{-1}*1}\\
0\to K^*_G(X)\ot K_*(A)\uuto^\cong_{\g\ot 1\qquad \fbox{3}}
\dto_\cong^{\f_\C\ot 1\qquad \qquad \fbox{5}}
\rto &K^*_G(X;A)\uto_\cong^\g \dto_{\f_A}& \\
0\to K^*_G(V)\ot K_*(A)\ddto_\cong^{\g\ot 1\qquad \fbox{4}}
\rto &K^*_G(V;A)\dto^\cong_\g & \\
&K_*((C_0(V)\ot A)\skr G)\dto^\cong_{\s_V}&\\
0\to K_*(C^V_G)\ot K_*(A)\rto^{\qquad \a_V}
&K_*(C^V_G\ot A)\rto^{\b_V\qquad }
\dto^{\cong\qquad \qquad \fbox{9}}_{\rm Bott}
&\Tor (K_*(C^V_G), K_*(A))\to 0 \dto_{\bigcap\quad}\\
&\bigoplus\limits_{q\in\Z_2}K_q(C^V_G\ot B)
\rto^{\iota_*}
& \bigoplus\limits_{q\in\Z_2}K_q(C^V_G\ot C),
&
\enddiagram
 $$
where the rows are exact: they present the main result of
\cite{SchoII}.  Let us consider the following diagram:
 $$
\diagram
&\bigoplus\limits_{q\in\Z_2}K_q(C^X_G\ot B)
\rto^{\iota_*}
\ddto^{\g\f_B\g^{-1}\qquad \fbox{10}}_{\fbox{6}\qquad }
& \bigoplus\limits_{q\in\Z_2}K_q(C^X_G\ot C)
\ddto^{\g\f_C\g^{-1}\quad\qquad  \fbox{7}}
& \Tor (K_*(C^X_G), K_*(A))\to 0
\lto_{\supset\qquad }\ddto^{\g\f_\C\g^{-1}*1}
\\
K^*_G(X;A)\urto_\cong^{{\rm (Bott)}\:\s_X\:\g\quad}
\drto_{{\rm (Bott)}\:\s_V\:\g\:\f_A\quad}&&& \\
&\bigoplus\limits_{q\in\Z_2}K_q(C^V_G\ot B)
\rto^{\iota_*}
& \bigoplus\limits_{q\in\Z_2}K_q(C^V_G\ot C)
& \Tor (K_*(C^V_G), K_*(A))\to 0. \lto_{\supset\qquad }
\enddiagram
 $$
Let us suppose that all squares
$\fbox{1},\dots,\fbox{10}$
are commutative. We obtain the diagram
\begin{equation}\label{diagrosnovtom}
\diagram
0\to K^*_G(X)\ot K_*(A)\rto
\dto_\cong^{\f_\C\ot 1\qquad \fbox{11}}
& K^*_G(X;A)\rto\dto^{\f_A\qquad \fbox{12}}
&\Tor(K_*(C^X_G),K_*(A))\to 0\dto_\cong^{\g\f_\C\g^{-1}*1}
\\
0\to K^*_G(V)\ot K_*(A)\rto
& K^*_G(V;A)\rto
&\Tor(K_*(C^V_G),K_*(A))\to 0.
\enddiagram
\end{equation}
Its rows are exact by the commutativity of
$\fbox{3}$ and $\fbox{4}$, and
 $\fbox{11}=\fbox{5}$. If $a\in K^*_G(X;A),$ then
 $$
 \begin{array}{ll}
(\g\f_\C\g^{-1}*1)\: (\b_X \circ\s_X\circ \g (A))=
&\qquad \fbox{7}\:\fbox{8}\\
\qquad =\g\f_C\g^{-1} \circ \iota_*\circ ({\rm Bott})\circ \s_X \circ
\g(a)=
 &\qquad \fbox{6}\:\fbox{10}\\
=\iota_*\circ ({\rm Bott})\circ \s_V \circ \g\circ \f_A(a)=
 &\qquad \fbox{7}\:\fbox{9}\\
=(\b_V\circ\s_v\circ\g)\f_A(a).
 \end{array}
 $$
Hence, (\ref{diagrosnovtom}) is a commutative diagram with exact rows.
Let us show that the squares $\fbox{1},\dots,\fbox{10}$ are
commutative.
 \begin{itemize}
  \item $\fbox{1}$ commutes by the standard algebraic argument;

  \item $\fbox{5}$ commutes by the associativity of
$\o$: rows and columns are induced by $\ot$;

  \item $ \fbox{2}$ commutes if we can prove the commutativity of
$\fbox{3}$ or $\fbox{4}$ for any unital $A$. Indeed,
$C$ is not unital, but we can consider the unitalization $C^+$ of it
and the diagram:
 \end{itemize}
 $$
\diagram
&\quad\qquad \qquad K_*(C^X_G)\ot K_*(C)\rto\dlto
& K_*(C^X_G\ot C)\dlto \\
 K_*(C^X_G)\ot K_*(C^+)\rto&K_*(C^X_G\ot C^+)\qquad \uto<-1.5em>\\
&\quad \qquad \qquad K^G_*(C_0(X))\ot K_*(C)\dlto\rto
\uline<-1.5em>
& K^G_*(C_0(X)\ot C)\dlto\uuto\\
K^G_*(C_0(X))\ot K_*(C^+) \rto\uuto
&K^G_*(C_0(X)\ot C^+)\qquad\uuto<2.5em> &\\
\enddiagram
 $$

 \begin{itemize}
  \item[] Then the forward square is commutative by
$\fbox{3}$. The commutativity of the others is evident.
Hence, the back square is commutative and we get
$\fbox{2}.$

  \item $\fbox{8}$ and $ \fbox{9}$ commute by the construction
of the exact sequence from \cite{SchoII}.

  \item $ \fbox{6}$ commutes since $\g$ for $K_1$ is defined in
 \cite{Green,Julg} via the Bott isomorphism.

  \item $ \fbox{7}=\fbox{1}+\fbox{2}$.

 \end{itemize}

So, we have to verify the commutativity of
$\fbox{3}$ (for any unital algebra) and $\fbox{10}.$

Let $P$ be a projective $G$-$C_0(X)$-module, $M$ a projective
$A$-module. Then
 $$
\a(\g\ot 1)([P]\otl_\Z[\M])=
\a([P\otl_{L^1(G,C_0(X))}(C_0\skr G)]\otl_\Z[\M])
=[(P\otl_{L^1(G,C_0(X))}(C_0\skr G)\otl_\C\M],
 $$
and
 $$
\s\g\a([P]\otl_\Z[\M])=\s\g([P\otl_\C\M])=
\s[(P\otl_\C\M)\otl_{L^1(G,C_0(X)\ot A)}((C_0(X)\ot A)\skr G)]=
 $$
 $$
=[P\otl_\C\M\otl_{L^1(G,C_0(X)\ot_{max} A)}
\bigl\{((C_0(X)\skr G)\ot_{max} A)\bigr\}].
 $$
Let us define the following mapping of the obtained modules
 $$
r:p\ot m\ot (f\ot a)\mapsto (p\ot f)\ot ma.
 $$
To prove that it is well defined, we have to verify
the following equality
 $$
r((p\ot m)h\ot (f\ot a))=r((p\ot m)\ot h(f\ot a)),
 $$
where $h\in L^1(G,C_0(X))\ot A$.
Without loss of generality it can be assumed that
 $$
h(g)=\sum_{i=1}^k h^i(g)\ot a^i,
 $$
where $h^i\in L^1(G,C_0(X))$, $a_i\in A$. Then
 $$
r((p\ot m)h\ot (f\ot a))
=\sum_{i=1}^k \int_G g^{-1}(p\cdot h^i(g))\,dg\ot f\ot m\cdot a^i a.
 $$
On the other hand,
 $$
r((p\ot m)\ot h (f\ot a))
=\sum_{i=1}^k p\ot h^i f \ot m\cdot a^ia
=\sum_{i=1}^k \int_G g^{-1}(p\cdot h^i(g))\,dg\ot f\ot m\cdot a^i a.
 $$
Hence, the mapping is well defined for the algebraic tensor
products. Since $A$ is a unital algebra,
$r$ is an isomorphism of the algebraic tensor products.
Consider
 $$
t:(P\otl_{L^1(G,C_0(X))}(C_0\skr G))\otl_\C\M \ola{}
(P\otl_\C\M)\otl_{L^1(G,C_0(X)\ot A)}((C_0(X)\skr G)\ot A),
 $$
 $$
t(p\ot f\ot m):=(p\ot m)\ot (f\ot 1_A).
 $$
Then
 $$
t\circ r (p\ot m \ot (f\ot a))=(p\ot ma)\ot(f\ot 1_A)=
\left(\frac 1\l p\ot m \right)\:(\l\ot a)\otl_{L^1}(f\ot 1_A)=
 $$
 $$
=\int_G g^{-1} \left(\frac 1\l p\ot m \right)\:(\l\ot a)(g)\, dg
\otl_{L^1}(f\ot 1_A)=\left(\frac 1\l p\ot m \right)\ot (\l f\ot a)=
(p\ot m) \ot (f\ot a),
 $$
i.~e., $t\circ r=\Id$.
We obtained an algebraic isomorphism of modules over
$(C_0(X)\skr G)\o A$. Let us consider its behavior with the respect
to norms. Considering a dense subset, we can assume that the element
has the form
 $$
x=\sum_{i=1}^k (p_i\ot m_i)\ot f_i \ot a_i,
 $$
where
 $$
f_i\ot a_i\in L^1(G,C_0(X))\ot_{max}A\ss (C_0(X)\skr G)\ot_{max}A.
 $$
On this dense subset $r$ is the following isometry:
 $$
P\otl_\C \M\ot \{1\}\ola{}(P\ot \{1\})\otl_\C\M.
 $$
Hence, $r$ is an isomorphism of $(C_0(X)\skr G)\ot A$-modules and
 $$
\s\g\a=\a(\g\ot 1).
 $$

Let us verify the commutativity of $\fbox{10}.$ By Lemma
\ref{lem:kommut} it is sufficient to obtain the equality
$\g\iota_*=\iota_*\g$. We have:
\begin{eqnarray*}
\g\iota_*P&=&(P\otl_{A\ot\K\ot C_0(\R)}C)\otl_{L^1(G,C_0(X)\ot C)}
((C_0(X)\ot C)\skr G)=\\
&=&(P\otl_{A\ot\K\ot C_0(\R)}C)\otl_{L^1(G,C_0(X))\ot_{max} C}
((C_0(X)\skr G)\ot C)=\\
&=&(P\otl_{A\ot\K\ot C_0(\R)}C)\otl_{L^1(G,C_0(X))}
(C_0(X)\skr G),
\end{eqnarray*}
and
\begin{eqnarray*}
\iota_*\g P&=&\left(P\otl_{L^1(G,C_0(X)\ot B)}(C_0(X)\ot B)\skr G\right)
\otl_{C_0(X)\ot B}C_0(X)\ot C=\\
&=&\left(P\otl_{L^1(G,C_0(X))\ot B}(C_0(X)\skr G)\ot B\right)
\otl_{B} C=\\
&=&\left(P\otl_{L^1(G,C_0(X))}(C_0(X)\skr G)\right)
\otl_{A\ot\K\ot C_0(\R)} C,
\end{eqnarray*}
where as before $B:=A\ot\K\ot C_0(\R)$. Let us define
 $$
\bar r:\g\iota_* P\ola{}\iota_*\g P,\qquad
\bar r(p\ot c\ot h)=p\ot h\ot c.
 $$
Then $\bar r$ is well defined. Indeed, if
$\a\in A\ot\K\ot C_0(\R)$ then
 $$
\bar r(p\a\ot c\ot h)=p\a\ot h\ot c,
 $$
 $$
\bar r(p\ot \a c\ot h)=p\ot h \ot \a c=p\a \ot h \ot c,
 $$
and if $f\in L^1(G,C_0(X))$ then
 $$
\bar r ((p\ot c)f\ot h)=\bar r (pf\ot c\ot h)=pf\ot h\ot c,
 $$
 $$
\bar r (p\ot c\ot fh)=p\ot fh\ot c=pf\ot h\ot c.
 $$
Let $f\ot d\in (C_0(X)\skr G)\ot C$. Then
 $$
\bar r ((p\ot c\ot h)(f\ot d))=\bar r (p\ot cd\ot hf)=p\ot hf\ot cd,
 $$
 $$
(\bar r (p\ot c\ot h))(f\ot d)=(p\ot c\ot h)\,(f\ot d)=p\ot hf\ot cd.
 $$
Hence, $\bar r$defines on the algebraic tensor product
a well-defined homomorphism of modules over $(C_0(X)\skr G)\ot C$. On
the dense subset
 $$
(P\otl_{A\ot\K\ot C_0(\R)} C)\otl_{L^1(G,C_0(X))} L^1(G,C_0(X))\ss
(P\otl_{A\ot\K\ot C_0(\R)} C)\otl_{L^1(G,C_0(X))} C_0(X)\skr G
 $$
the homomorphism $\bar r$ is the following isometry
 $$
(P\otl_{A\ot\K\ot C_0(\R)} C)\ot \{1\}\to
(P\ot \{1\})\otl_{A\ot\K\ot C_0(\R)} C.
 $$
Hence,  $\bar r$ is an isomorphism of
modules over the algebra $(C_0(X)\skr G)\o C$ completed with
respect to $\ot_{\rm max}$; and $\g\iota_*=\iota_*\g.$

Now we apply Five-lemma to the diagram (\ref{diagrosnovtom}).
The proof of Theorem is completed.
\hfill $\Box$

 \begin{th}\label{isotomatriv}
If $X$ is a separable metrizable trivial $G$-space, then
$\varphi_A$ is an isomorphism.
 \end{th}

\Dok Let us consider the diagram:
 $$
 \begin{array}{cc}
\bK^*_G(X;A)\cong
&\bK^G_*(C(X)\ot A)\stackrel\g{\cong}
K_*((C(X)\ot A)\skr G)\cong
\\
\downarrow\:{\f_A}&\\
\bK^*_G(V;A)\cong
&\bK^G_*(C(V)\ot A)\stackrel\g{\cong}
K_*((C_0(V)\ot A)\skr G)\cong
 \end{array}
 $$
 $$
 \begin{array}{ccc}
\cong
K_*(C(X)\ot (A\skr G))\cong&
K^*(X;A\skr G)\\
& \cong\:\downarrow\:{\f_{A\skr G}}\\
\cong K_*(C_0(V)\ot (A\skr G))\cong&
K^*(V;A\skr G).
 \end{array}
 $$
To prove the commutativity of this diagram (which implies Theorem)
it is sufficient to demonstrate the following.
Let $\M$ be a projective $GG,C_0(V)\ot A$-module of finite type and
$E$ is a projective $GG,C(V)$--module, i.~e.
$C(V)$-module, because the action of $G$ on $C(V)$ is trivial.
Then
 $$
\g(\M\ot_{C(V)}E)=\g(\M)\ot_{C(V)}E.
 $$
Since these are ``the same modules'' by the definition of $\g$ in
\cite{Julg} (where it was denoted by $\Psi$), it is necessary to verify
first, that the module structure survives, and second, that the
constructions behaves properly with the respect to morphisms
$f:E\to F$. Indeed,
 $$
(m\ot e)\cdot_{\g(\M\ot_{C(V)}E)}
\left(\int_G a(g) U_g dg \ot \f\right)=
\left(\int_G g^{-1}[(m\ot e)\cdot (a(g)\ot \f)  dg \right)=
$$
$$
=
\left(\int_G g^{-1}[(m\cdot a(g))\ot (e\cdot \f) dg \right)=
$$
$$
=
\left(\int_G g^{-1}(m\cdot a(g)) dg \right)\ot (e\cdot \f)=
(m\ot e)\cdot_{\g(\M)\ot_{C(V)}E} \left(\int_G a(g)U_g dg\ot
\f \right).
 $$
The statement about morphisms is immediate consequence of the
triviality of the action of $G$ on $C(V)$.
\hfill $\Box$

Now, after the particular cases, we are able to prove the general
theorem about the Thom isomorphism.

 \begin{th}\label{isotomaobsh}
Let $X$ be a manifold, then $\f_A$ is an isomorphism.
 \end{th}

\Dok
First, let us prove the theorem for the trivial bundle $X\times V$,
where $V$ is a complex finite-dimensional $G$-space.

Let us denote by $1$ the trivial $1$-dimensional $G$-module, and the
projective space $P(V\oplus 1)$ is a compactification of $V$ and we have
the following natural homomorphism
 $$
j:\bK_G(V;A)\to \bK_G(P(V\oplus 1);A),\qquad
j:\bK_G(V\times X;A)\to \bK_G(P(V\oplus 1)\times X;A).
 $$

Let $X$ be compact. Let us consider an arbitrary element
$x\in \bK_G(P(V\oplus 1)\times X; A)$. Let us consider the analytical
index of the correspondent family of Dolbeault operators over
$P(V\oplus 1)$ with coefficients in $x$
(cf.~\cite{At-Period}), i.~e. an operator over $C^*$-algebra $C(X)\ot A$
(see Section~\ref{sec:famil}). This is an element of
$\bK_G(X;A)$. Taking the composition with $j$
(cf.~\cite[p.~123]{At-Period}) we get a family of mappings
$\a=\a_{X,A}:\bK_G(V\times X;A)\to \bK_G(X;A)$, having the following
properties:
 \begin{description}
  \item[(a1)] $\a$ is functorial with the respect to $G$-morphisms of $X$
and $A$;

  \item[(a2)] the following diagram commutes:
 $$
\diagram
\bK_G(V\times X;B)\ot \bK_G(X;A)\rto\dto^{\a_{X,B}\ot 1} &
\bK_G(V\times X; B\ot A)\dto^{\a_{X,B\ot A}}\\
\bK_G(X;B)\ot\bK_G(X;A)\rto&\bK_G(X;B\ot A),
\enddiagram
 $$
in particular, $\a$ is a morphism of $K_G(X)$-modules;

  \item[(a3)] $\a_{pt,\C}(\la_V^*)=1\in R(G)$.

 \end{description}

In fact, transferring the space of parameters into coefficients:
 $$
\bK_G(P(V\oplus 1)\times X;B)\cong \bK_G(P(V\oplus 1);C(X)\ot B),
 $$
we reduce the situation to the case $X=pt$, and (a2) takes form
\begin{equation}\label{a21}
\diagram
\bK_G(V;B)\ot \bK_G(pt;D)\rto\dto^{\a_{pt,B}\ot 1} &
\bK_G(V; B\ot D)\dto^{\a_{pt,B\ot D}}\\
\bK_G(pt;B)\ot\bK_G(pt;D)\rto&\bK_G(pt;B\ot D),
\enddiagram
\end{equation}
and the functoriality of (a1) means functoriality with the respect to the
algebra of coefficients. In a usual way with the help of commutative
diagram
 $$
\diagram
0\to\bK_G(V;A)\rto&\bK_G(V;A^+)\rto\dto^{\a}
&\bK_G(V;\C)\dto^\a\\
0\to\bK_G(pt;A)\rto &\bK_G(pt;A^+)\rto&\bK_G(pt;\C)
\enddiagram
 $$
$\a$ can be extended to non-unital (in particular, non-compact) case.
The diagram (\ref{a21}) is still commutative. We will need the following
two particular cases of it.
For $B=C_0(X)\ot A$, $D=C_0(V)$ we get for
$x\in \bK_G(V\times X;A),\quad y\in K_G(V)$
\begin{equation}\label{a4}
\a_{X\times V, A}(x\ot y)=\a_{X,A}(x)\ot y\in
\bK_G(X\times V; A).
\end{equation}
For $D=C_0(X)\ot A$, $B=C_0(W)$ we get for
$y'\in \bK_G(W\times X;A),\quad x'\in K_G(V)$
\begin{equation}\label{a5}
\a_{W\times X, A}(x'\ot y')=\a_{pt,\C}(x')\ot y'\in
\bK_G(W\times X; A).
\end{equation}
Let $x\in \bK_G(X;A)$, then by (\ref{a5}) with $W=0$ and (a2)
\begin{equation}\label{a6}
\a(\la_V^* x)=\a(\la_V^*) x=x.
 \end{equation}
Let $y\in \bK_G(V\times X;A)=\bK_G(V; C_0(X)\ot A)$, then by (\ref{a4})
and (\ref{a5}) with $W=V$
\begin{equation}\label{a7}
\a(y)\ot\la_V^*=\a(y\ot \la_V^*)=
\a(\la_V^* \ot\wt  y)=\a(\la_V^*) \wt y=\wt y\in \bK_G(X\times V;A),
\end{equation}
where by $\wt y\in \bK_G(X\times V;A)$ we denote the element, obtained
from $y$ under the mapping
$X\times V \to V \times  X$, $(x,v)\mapsto (-v,x)$
(such that $V\times X \times V\to V\times X \times  V$,
$(u,x,v)\mapsto (-v,x,u)$ is homotopic to the identity). Let us
apply to the both parts of (\ref{a7}) the isomorphism
$\bK_G(X\times V;A)\to\bK_G(V\times X;A)$. We obtain that
$\f\a$ is an isomorphism.  But by (\ref{a6}) $\a\f=\Id$, hence,
$\a$ is the two-sides inverse. And the unknown automorphism is the
identity.

The pass to the case of general complex $H$-bundle $E\to Y$ we
make similarly to \cite[p.~124]{At-Period}. Namely, let us take
$G=U(n)\times H$ and the principal
$G$-bundle $X$, with which $E$ is associated, $Y=X/U(n)$,
in such a way that
 $$
\bK_G(X;A)\cong \bK_H(Y;A),\qquad
\bK_G(V\times X;A)\cong \bK_H((V\times X)/U(n);A) \cong \bK_H(E;A).
 $$
Since the identifications commute with multiplication by
$\la_V^*$, we get the desired result.
(Other way is to enclose the bundle into a trivial one and use
the transitivity of the Thom homomorphism).
\qquad $\Box$

We shall pass now to some further constructions, connected with
the Thom homomorphism. They will be necessary for the definition of
topological index.

Let $X$ and $Y$ be smooth $G$-manifolds, $i: X\to Y$
the equivariant enclosure, $Y$ is equipped with a $G$-invariant
Riemannian metric, $(TX,p_T,X)$ is the tangent bundle of $X,$
$(N,p_N,X)$ is the normal bundle for $i.$ Let us choose a
function $\e: X\to (0,\infty)$ such that the map of $N$ to itself
 $$
n\mapsto \e\,\frac n{1+|n|}
$$
is $G$-equivariant and determines a
$G$-diffeomorphism $\F: N\to W$ on an open tubular neighborhood
$W\supset X$ in $Y.$ The enclosure $i: X\to Y$ is decomposed in a
composition of two enclosures $i_1: X\to W$ and $i_2: W\to Y.$ Passing to
differentials we obtain
 $$
TX\ola{di_1}TW\ola{di_2}TY,\qquad d\F:TN\to TW.
 $$
 \begin{lem}
{\em \cite[p. 112]{Fried}\/}
The manifold $TN$ can be identified with $p^*_T (N\oplus N)$
with the help of a $G$-dif\-fe\-o\-mor\-phism $\psi$ such that the
following diagram is commutative
 $$
\diagram
p^*_T(N\oplus N)\dto && TN \llto^\psi\dto\\
TX\drto^{p_T} &  & N\dlto_{p_N}\\
&X.&
\enddiagram
 $$
 \end{lem}

\Dok
The manifold $TN$ as the vector bundle over $N$
can be identified with $p^*_N(TX)\oplus p^*_N(N).$
A point of the total space $TN$ is a pair of the form
$(n_1, t + n_2),$ where both vectors are from the fiber
over the point $x\in X.$ Similarly, we represent elements
$p^*_T(N\oplus N)$ as pairs of the form $(t,n_1 + n_2).$
Let us define $\psi$ by the equality
$\psi (n_1, t+n_2)= (t, n_1+n_2).$ \hfill$\Box$

With the help of the relation $i\cdot (n_1,n_2)=(-n_2,n_1)$,
we can equip
$$
p^*_T(N\oplus N)=p^*_T(N)\oplus p^*_T(N)
$$
with a structure of a complex\label{ccc} manifold.
Then we can consider the Thom homomorphism
 $$
\f:\bK_G(TX;A)\to \bK_G(p^*_T(N\oplus N),A).
 $$
Since $TW$ is an open $G$-stable subset of $TY$
and $di_2: TW\to TY$ is an enclosure, by the construction
\ref{conext}, there is the homomorphism
 $
(di_2)_*:K_G(TW;A)\to K_G(TY;A).
 $

\dfi\label{oprgys}
Let $i: X\to Y$ be an enclosure. The {\em Gysin homomorphism\/}
\label{ind:Gysinhomomorphism} is the mapping
 $$
i_!:\bK_G(TX;A)\to \bK_G(TY;A),\qquad i_!=(di_2)_*(d\F^{-1})^*\psi^*\f.
 $$
In other words, it is obtained by passage to $K$-groups in the
upper part of the diagram
 $$
\diagram
p^*_T(N\oplus N)\dto_{q_T} &&TN\llto_\psi\rto^{d\F}\dto
& TW \rto^{di_2}\ddto& TY\ddto \\
TX \drto^{p_T} & & N\dlto_{p_N} \drto^{\F} & & \\
 & X\rrto^{i_1} & & W \rto^{i_2} & Y.
\enddiagram
 $$

Another choice of metric and neighborhood $W$
induces the homotopic map and (by the item 3
of Theorem~\ref{th:prophys} below) the same homomorphism.

 \begin{th}\label{th:prophys}
{\em (Свойства гомоморфизма Гизина)\/}
\noindent
{\rm\bf 1.} $i_!$ is a homomorphism of $R(G)$-modules.

\smallskip\noindent
{\rm\bf 2.} Let $i: X\to Y$ and $j: Y\to Z$ be two $G$-inclusions, then
 $
(j\circ i)_!=j_!\circ i_!.
 $

\smallskip\noindent
{\rm\bf 3.} Let enclosures $i_1: X\to Y$ and $i_2: X\to Y$
are $G$-homotopic in the class of enclosures. Then
 $
(i_1)_!=(i_2)_!.
 $

\smallskip\noindent
{\rm\bf 4.} Let $i_!: X\to Y$ be a $G$-diffeomorphism, then
 $
i_!=(di^{-1})^*.
 $

\smallskip\noindent
{\rm\bf 5.}
An enclosure $i: X\to Y$ can be represented as a
compositions of enclosures $X$ in $N$ (as the zero section $s_0: x\to N$)
and $N\to Y$ by $i_2\circ \F: N\to Y.$ Then
 $
i_!=(i_2\circ \F)_! (s_0)_!.
 $

\smallskip\noindent
{\rm\bf 6.} Consider the complex bundle $p^*_T (N\o\C)$ over $TX.$
Let us form the complex $\La(p^*_T(N\o\C),0):$
$$
0\to\La^0(p^*_T(N\o\C))\ola{0}\dots
\ola{0} \La^k(p^*_T(N\o\C))\to 0
 $$
with the noncompact support. If $a\in \bK_G(TX;A),$ then the complex
$
a\o \La(p^*_T(N\o\C),0)
$
has the compact support and determines an
element of $\bK_G (TX; A).$ Then
$
(di)^* i_! (a)=a\cdot \La(p^*_T(N\o\C),0).
$

\smallskip\noindent
{\rm\bf 7.} For $x\in \bK_G(TX;A)$ and $y\in K_G(TY)$
$
i_!(x(di)^*y)=i_!(x)\cdot y.
$
 \end{th}

\Dok
1. By the definition of $i_!.$

2. To simplify the argument, let us identify the tubular
neighborhood with the normal bundle. Then $(j\c i)_!$ is the
composition
 $$
\bK_G(TX;A)\ola{\f} \bK_G(TN\oplus TN'_1;A)\to \bK_G(TZ;A),
 $$
where $N'$ is the normal bundle of $Y$ in $Z,$ $N'_1=N'|_X,$ and
$TN\oplus TN'_1$ is considered in the same way as at p.~\pageref {ccc}
as a complex bundle over $TX.$ On the other hand,
$j_!\c i_!$ represents the composition
$$
\bK_G(TX;A)\ola{\f} \bK_G(TN;A) \to \bK_G(TY;A)
\ola{\f} $$
$$ \ola{\f} \bK_G(TN';A) \to \bK_G(TZ;A).
$$
By properties of $\f$, the following diagram is commutative
$$
\begin{array}{ccccc}
 \bK_G(TX;A)&\ola{\f}& \bK_G(TN;A)& \to& \bK_G(TY;A)\\
&\searrow{\scriptstyle\f} &\downarrow{\scriptstyle\f}
& & \downarrow {\scriptstyle\f}\\
&& \bK_G(TN\oplus TN'_1;A) & \ola{\f} & \bK_G(TN';A) \\
&&\downarrow &&\downarrow \\
& & \bK_G(TZ;A) &= & \bK_G(TZ;A),
\end{array}
$$
This completes the proof of item 2.

3. The homotopy of enclosures does not influence on $q_T,$ but only on
the further maps in the definition of the Gysin homomorphism. But for
these maps the assertion follows from the homotopy invariance
of $K$-theory.

4. In this case $N=X,\,W=Y,\,\F=i,\,i_2=\Id_Y,$ and
the formula is obvious.

5. By 2.

6. By definition,
 $$
(di)^*i_!=(di_1)^*(di_2)^*(di_2)_*(d\f^{-1})^*\psi^*\f^* =
(\psi\c d\F^{-1}\c di_1)^*\c \f,
 $$
where $i_1: X\to W,\, i_2: W\to Y.$ Let
$(n_1,t+n_2)\in TN=p^*_N(TX)\oplus p^*_N(N)$,
where $n_1$ is the shift under the exponential mapping,
$t+n_2$ a tangent vector to $W.$
 If $d\F(n_1,t+n_2)$ is in $TX, then $ $n_1=n_2=0.$ Hence,
 $$
d\F^{-1}di_1(t)=(0,t+0),\quad
\psi\c d\F^{-1}\c di_1 (t)=(t,0+0).
 $$
Therefore,
 $
\psi\c d\F^{-1}\c di_1 : TX \to p^*_T(N\oplus N)
 $
is the enclosure of the zero section. Since,
 $
\f(a)=a\cdot \La(q^*_T p^*_T (N\o \C), s_{p^*_T(N\o\C)}),
 $
it follows that
 $
(di)^*i_!(a)=a\cdot \La(p^*_T(N\o\C),0).
 $

7. The mapping $di_1\c q_T \c\psi\c d\F^{-1}:TW\to TW$
is homotopic to the identical mapping. Hence,
$$
\begin{array}{l}
 i_!(x\cdot (di)^*y)=(di_2)_*(d\F^{-1})^*\psi^*\f(x\cdot (di)^*y)=\\
 \:=(di_2)_*(d\F^{-1})^*\psi^*[(q^*_T(x)\la_{p^*_T(N\o\C)})
(q^*_T(di)^*y)]=\\
 \:=(di_2)_*[(d\F^{-1})^*\psi^*(q^*_T(x)\la_{p^*_T(N\o\C)})
\underbrace{(d\F^{-1})^*\psi^*q^*_T(di_1)^*}_{\Id}(di_2)^*y]=\\
 \!=[(di_2)_*(d\F^{-1})^*\psi^*(q^*_T(x)\la_{p^*_T(N\o\C)})]\,
[(di_2)_*(di_2)^*y]=i_!(x)\cdot y. \quad\Box
         \end{array}
$$

 \begin{th}\label{th:prophys2}
{\rm\bf 1.} Let $V$ be a
$G$-$\R$-vector space and $X=pt$ a trivial
$G$-manifold. Hence, $TX=pt$ and $TV=V\o\C.$
Consider the enclosure of zero $i: X\to V$. Then
the mapping
 $$
i_!:\bK^G(A)=\bK_G(TX;A)\to \bK_G(TV;A)=\bK_G(V\o\C;A)
 $$
coincides with the Thom homomorphism $\f^{V\o\C}_A.$

\smallskip\noindent
{\rm\bf 2.}
Let $V_1$ and $V_2$ be $G$-$\R$-spaces, $i: X\to V_1$ an
enclosure. Let us define the enclosure
$k: X\to V_1\oplus V_2$ by the formula $k(x) =i(x) + 0.$ Then
the following diagram is commutative
 $$
\diagram
& \bK_G(TV_1;A)\ddto^\f \\
\bK_G(TX;A)\urto^{i_!}\drto_{k_!}& \\
& \bK_G(T(V_1\oplus V_2),A),
\enddiagram
 $$
where $\f$ the Thom homomorphism of the complex bundle
 $$
T(V_1\oplus V_2)=V_1\o\C \oplus V_2\o\C
\ola{}TV_1=V_1\o\C.
$$
 \end{th}

\Dok
1. The assertion follows from the definition of $i_!.$
More precisely, $X=0\in V,\, N=V$. Also, $W$
can be chosen equal to the interior $D_1$ of
the ball of radius 1 in $V$ with respect to an
invariant metric.
In this case the diagram from the definition of the Gysin
homomorphism \ref{oprgys} takes the following form
 $$
\diagram
V\o\C\dto\rrto^\Psi &&V\o\C\rto^{d\F}\dto
& D_1\o\C \rto^{di_2}\ddto& V\o\C\ddto \\
TX=0 \drto & & V\dlto \drto^{\F} & & \\
 & X=0\rrto^{i_1} & & D_1 \rto^{i_2} & V.
\enddiagram
 $$
In our case $\Psi=\Id$ and $di_2\c d\F$ is homotopic to $\Id,$
since it has the form
 $
v\o z \mapsto  (v\o z)/({1+|v\o z|}).
 $
Hence, $i_!=\f.$

2. Let the enclosure $i: X\to V_1$ has the normal bundle $N$ and
the tubular neighborhood $W.$ Then $N\oplus V_2$ is a normal
bundle for the enclosure
$k$ with the tubular neighborhood $W\oplus D (V_2),$ where
$D(V_2)$ is a ball. If $a\in \bK_G (TX; A),$ then
 $$
k_!(a)=(di_2\oplus 1)_*(d\F^{-1}\oplus 1)^*(\psi\oplus 1)^*
\f_A^{N\oplus N \oplus V_2 \oplus V_2}(a).
 $$
By item 2 of theorem \ref{svolamb},
 $
\f_A^{N\oplus N \oplus V_2 \oplus V_2}=
\f_A^{N\oplus N}\cdot\f_\C^{V_2 \oplus V_2}.
 $
Since $a=a\cdot \underline {\C},$ where $\underline {\C} $
is the trivial line bundle,
 $$
k_!(a)=(di_2)_*(d\F^{-1})^*\Psi^*\f^{N\oplus N}_A (a)\cdot
\f_\C^{V_2 \oplus V_2}(\underline{\C})=
 $$
 $$
=i_!(a)\cdot \la_{T(V_1\oplus V_2)}=\f (i_!(a)).\qquad\Box
 $$

 \section{Analytical index}

let us remind how a
pseudo differential operator $\chi(\s)$ over $A$
can be defined starting from a symbol $\s$ in the case $G=e$
(see, e.g.~\cite{SolTro}).

Suppose $X$ is a compact closed smooth manifold,
$\pi: T^*X\to X$ is the projection of the cotangent bundle, $E_1$
and $E_2$ are smooth $A$-bundles over $X$. Let $\{U_j\}$ be
a trivializing cover of $X$ for $E_1$ and $E_2$, $\{\f_j\}$ be
the subordinate partition of unity, $\{\psi_j\}$ a collection of
functions such that $\psi_j|_{\supp \f_j}\equiv 1$.
A symbol of order $m$ is a morphism of
$A$-bundles $\s:\pi^*E_1\to\pi^*E_2$.

Let $u_j\in\G^\infty_0(U_j,E_1)$ be a smooth section tending to
$0$ at infinity.  Let us assume
 $$
 (\chi_j\s)(u_j)(x):=\frac
1{(2\pi)^n}\int e^{i\<x,\xi\>}\s(x,\xi)\int
 e^{-i\<y,\xi\>}u_j(y)\,dy\,d\xi.
  $$
For $u\in\G^\infty (X,E_1)$ let us define $\chi\s$ by the formula
 $$
(\chi\s)(u)(x)=\sum_j \psi_j(x)(\chi_j\s)(\f_j u)(x).
  $$

This is the section of the symbol-map  (see, e.~g.
\cite[\S 2.1]{SolTro}).
We obtain a bounded $A$-homomorphism
admitting an adjoint in Sobolev spaces:
 $$
\chi\s:H^s(E_1)\cong l_2(\P_1)\to l_2(\P_2)\cong H^{s-m}(E_2).
  $$

Suppose $G$ is an arbitrary compact Lie group,
$X$ is a compact $G$-manifold,  $E_1$ and $E_2$ are
 $GGA$-bundles over $X$, the symbol $\s:\pi^* E_1\to \pi^* E_2$ is a
 $GGA$-morphism.
The action of $G$ on Sobolev spaces is naturally defined by the
formula $g(u)(x)=g(u(g^{-1}x))$.
Without loss of generality, we can assume the
Sobolev inner products
on $H^s (E_i)\cong l_2 (\P_i)$ to be $GGA$-products
(see Remark~\ref{usrP}).

More precisely, let us choose a
$C^*$-Hermitian metrics in bundles to be $G$-invariant, we obtain
that under the action of $G$ on Sobolev spaces we get each time
admissible products (see, e.~g.
\cite[\S 2.1]{SolTro}). Hence, this $\Psi$DO admits an adjoint with the
respect to the averaged Sobolev product.

Suppose  $\s$ is an elliptic symbol, then (see, e.~g.
\cite[\S 2.1]{SolTro})
 \begin{equation}\label{3.2}
 (\chi\s)\,(\chi\tau)=1+K_1,\qquad (\chi\tau)\,(\chi\s)=1+K_2
 \end{equation}
for some symbol (parametrix) $\tau$, where $K_1$ and $K_2$ are
$A$-compact operators.
The action of group $G$ on $\Psi$DO is given by $g(P)u=g(P(g^{-1}(u)))$.

 Let $g\in G$ be fixed.
Let us denote by $g\{\:\} $ the action of $g$ on
sections and by $g\cdot$ the action on the elements of
total space of a bundle.
By~\cite[\S 2.1]{SolTro} $\OP_{r-1}$ and $CZ_r$,
hence, $\Int_r$ are invariant under diffeomorphisms,
in particular, under $g$.
Hence $\Int_r$ has this invariance property. Moreover
(see, e.~g. \cite[\S 2.1]{SolTro}),
 $$
\s_m(g(\chi\s))(x,\xi)=\s_m(g\cdot\chi\s\cdot g)(x,\xi)=
$$
$$
=\s_m(\chi\s)(gx,g\xi)=\s(gx,g\xi)=\s(x,\xi),
 $$
since $g$ acts on covectors by
 $$
g\xi=g_*\xi=(D(g^{-1}))^*\xi,
 $$
where the tangential mapping (derivative) of an arbitrary
smooth map $\f$ we denote by $D(\f)$. Hence
 $$
g(\chi\s)=g\{\chi\s\}g^{-1}=\chi'\s,
  $$
i.~e., $g (\chi\s) $ is
another $\Psi$DO with the same symbol $\s$, where $\chi'$ is some other
 section of the symbol-map.

  \begin{lem}
Let us fix $D=\chi\s$. Then the map $G\to \Int_m (E,F)$,
given by the formula
 $$
g\mapsto g(D),\qquad g(D)(u)=g\{Dg^{-1}\{u\}\},
  $$
is continuous, i.~e., for any $s$ this map is
continuous as a map from the group $G$ to the space of operators
 $H^s\to H^ {s-m}$ equipped with the uniform topology.
  \end{lem}

\Dok (cf.~\cite{AZ1})
Let $a\in \gg$, where \gg\ is the Lie algebra
of group $G$. The action of $a$ defines
a vector field  $\vec a_x$. The corresponding
differential operator $a_E$ acts on
sections of the bundle $E$:
 $$
 a_E(u)=\left\<\vec a_x,\frac {\vec\partial}{\partial x_i}\right\>(u).
 $$
The symbol of this operator is equal to $\s_{a,E}=\<\vec a,\xi\>\Id_E$.
Similarly for $F$. Then $\s\s_{a,E}=\s_{a,F}\s$
(see, e.~g. \cite[\S 2.1]{SolTro}),
$$
 D\,(\chi\s_{a,E})-(\chi\s_{a,F})\,D\in \Int_m(E,F),
$$
$$
\|D\,(\chi\s_{a,E})-(\chi\s_{a,F})\,D\|^m_s\le\|D\|(\|a_E\|+\|a_F\|),
  $$
where $a_E=\chi\s_{a,E}$. If $a$ is in a bounded neighborhood of
 $0\in \gg$, there exist numbers $c_s$  such that
 $$
\|D\,a_E-a_F\,D\|^m_s < c_s,
  $$
where $\|.\|^m_s$ is the norm in $\Int_m$, i.~e., the uniform norm
in the space of operators $H^s\to H^{s-m}$.

For $g_t=\exp ta$ and a smooth section $u$ let us assume
$f (t) :=g_t (D) u$.
Then
$$
f_t(t)=\exp (ta_F)D\exp (-ta_E)u.
$$
Indeed, it is sufficient to prove this equality
in one chart for analytical functions, because
the analytical basis was constructed in an explicit form
in~\cite[p. 854]{MF}.
For an analytical section $u(x)$ we have to prove that
\begin{equation}\label{3.3}
(e^{ta_E}u)(x)=u(\f_x(t)),
\end{equation}
where $\f_x (t) $ is the orbit of the action of the one-parameter
subgroup, generated by $a$. For this purpose let us
consider the function $F (t) =u (\f_x (t)) $. It is analytical
for $|t|\le 1$.  Then for
 $|t|\le 1$ the series
 $$
F(t)=\sum_{n=0}^\infty \frac {F^{(n)}(0)}{n!}t^n
 $$
is convergent. Since $\f_x (t)$ is an integral curve of the
vector field $\vec a_x$,
 $$
F'(t)=\frac {du(\f_x(t))}{dt}=\frac {d\f_x(t)}{dt}u=\vec a_{\f_x(t)}u=
(a_E u)(\f_x(t)),
 $$
i.~e.,  $F'(t)$  plays the role of
$F(t)$ for the section $a_E u$. By induction
 $$
F^{(n)}(t)=(a^n_E u)(\f_x(t)).
 $$
Therefore
 $$
F^{(n)}(0)=(a^n_E u)(\f_x(0))=(a^n_E u)(x).
 $$
Hence
 $$
F(t)=\sum_{n=0}^\infty \frac {(a^n_E u)(x)}{n!}t^n=
\sum_{n=0}^\infty \frac {(ta_E)^n u}{n!}.(x)=(e^{ta_E}u)(x),
 $$
and we have (\ref{3.3}).
Further, since the action of $G$ on $H^s$ is an isometry,
 $$
\left\|\frac {df_t}{dt}\right\|_{m-s}=
\|\exp ta_F\,(D\,a_E-a_F\,D)\,\exp (-ta_E)u\|\le c_s \|u\|_s.
 $$
This implies continuity.
\hfill$\Box$

Let $\chi^1\s=D_1$ and $\chi^2\s=D_2$ be two $\Psi$DO with the same
symbol $\s$. Then $D_1-D_2\in\K (l_2(\P_1,\P_2))$.
If $\s$ is a $G$-equivariant symbol, then we have shown that
$g (D) =D'=\chi'\s$. Hence the orbit $G(D)$ is in the
closed linear manifold
 $$
D+\K(l_2(\P_1),l_2(\P_2)).
 $$
By the previous lemma there exists the averaged operator ${\rm Av}\,D$
and it lies in this linear manifold.
Thus for an elliptic $G$-symbol the operator
${\rm Av}\,D$ is an $A$-Fredholm one.
 Really, from (\ref{3.2}) it follows
 $$
{\rm Av}\,(D)\:{\rm Av}\,(\chi\tau)=(\chi\s+K'_1)(\chi\tau+K''_1)=
 $$
$$
=(\chi\s)\,(\chi\tau)+K'_1(\chi\tau)+ (\chi\s)\,K''_1 + K'_1 K''_1=
 1+\widetilde K_1,
 $$
 $$
 {\rm Av}\,(\chi\tau)\:{\rm Av}\,(D)=(\chi\tau+K''_1)(\chi\s+K'_1)=
 $$
$$
=(\chi\tau)\,(\chi\s)+K''_1\,(\chi\s)+(\chi\tau)\,K'_1 + K''_1\, K'_1=
1+\widetilde K_2,
$$
where $K'_1,\,K''_1,\,\widetilde K_1,\,\widetilde K_2$ are
$A$-compact operators. It was proved that each $\Psi$DO over
unital C*-algebra admits an adjoint. Hence
each element of the orbit $G(D)$ has an adjoint. Since the operators
admitting an adjoint form a Banach space (algebra), the averaged
operator ${\rm Av}\,D$ has an adjoint.
It remains to apply our theory of Fredholm operators.
Thus the equivariant analytical index
$$
\ai \s=\index {\rm Av}\,D\in
K^G(A)=\bK^0_G({\rm pt};A).
$$
is defined. It is clear, that
all argument about homotopies (see, e.~g. \cite{SolTro})
is valid in the equivariant
case. Hence in this case
the analytical index defines a homomorphism of Abelian
groups
 $$
\ai^X_G:\bK_G(TX;A)\to K^G(A)=\bK^0_G({\rm pt};A).
 $$

 \section{The axiomatic approach}

Let us start from the definition of the topological index.
Let $G$ be a compact Lie group and $X$ be a compact $G$-manifold.
From \cite{PalImb} it follows that
there exists a representation of $G$ in orthogonal group $O V)$ of some
real finite-dimensional space $V$ and $G$-enclosure $i: X\to V$.
Thus, the Gysin homomorphism (see~\ref{subs:Thom}):
 $$
i_! : \bK_G(TX;A)\to \bK_G(TV;A) = \bK_G(V\o\C;A)
 $$
is defined .
Since $TV=V\o\C$ is a complex vector space,
we have the following Thom isomorphism (see~\ref{subs:Thom}):
$$
\f:K^G_0(A)=\bK^0_G(pt;A)\to \bK_G(TV;A).
$$

\dfi\label{dfi:topind}
The {\em topological index\/}
is the following map:
 $$
\ti^X_G:\bK_G(TX;A)\to K^G(A),\qquad \ti^X_G:=\f^{-1}\circ i_!.
 $$

 \begin{th}\label{th:ax2} \begin{enumerate}
  \item The index $\ti^X_G$ does not depend on the choice
        of $V$, enclosure
        $i: X\to V$, and representation $G\to O (V)$.

  \item The index $\ti^X_G$ is a $R(G)$-homomorphism.

  \item If $X=pt$, then the map
 $$
\ti^X_G:K^G(A)=\bK_G(TX;A)\to K^G(A)
 $$
coincides with $\Id_{K^G(A)}$.

  \item Suppose $X$ and $Y$ are compact $G$-manifolds, $i: X\to Y$ is
	a $G$-enclosure. Then the diagram
 $$
\diagram
\bK_G(TX;A)\rrto^{i_!}\drto_{\ti^X_G\:\:} && \bK_G(TY;A)
\dlto^{\ti^Y_G}\\
&K^G(A).&
\enddiagram
 $$
commutes.
 \end{enumerate}
\end{th}

\Dok
1). Let us consider the enclosures
 $$
i_1:X\to V_1,\qquad i_2:X\to V_2.
 $$
Denote by $j=i_1+i_2$ the induced enclosure $j: X\to V_1\oplus V_2$.
It is sufficient to show that the topological index, which comes from
$i_1$, coincides with the index, which comes from $j$. Let us define
a homotopy of $G$-enclosures by the formula
 $$
j_s(x)=i_1(x)+s\cdot i_2(x)\::\:X\to V_1\oplus V_2,\quad 0\le s\le 1.
 $$
Then by Theorems~\ref{th:prophys}.3 and \ref{th:prophys2}.1 the
indexes for $j$ and $j_0$ coincide. Let us show now
that $j_0=i_1 + 0$ and $i_1$ define the same topological
indexes. For this purpose consider the diagram
 $$
\diagram
&\bK_G(TX;A)\dlto_{(i_1)_!}\drto^{(j_0)_!}&\\
\bK_G(TV_1;A)\rrto^{\f_2} && \bK_G(T(V_1\oplus V_2);A), \\
&K^G(A)\ulto^{\f_1}\urto_{\f_3}&
\enddiagram
 $$
where $\f_i$ are the corresponding Thom homomorphisms.
The upper triangle is commutative by Theorem \ref{th:prophys2}.2,
and the lower is commutative by Proposition \ref{prop:thom6}. Hence
$\f_1^{-1}\c(i_1)_!=\f_3^{-1}\c (j_0)_1$ as desired.

2). This statement follows from \ref{prop:thom3} and \ref{th:prophys}.1.

3). This follows immediately from the definition of the index and
	from\ref{th:prophys2}.1

4). Let us consider the diagram
 $$
\diagram
X\rrto^i\drto_{j\circ i}&&Y\dlto^j\\
&V.&
\enddiagram
 $$
Let us apply \ref{th:prophys}.2. We have
the following commutative diagram
 $$
\diagram
\bK_G(TX;A)\rrto^{i_!}\drto_{(j\circ i)_!\:}&&\bK_G(TY;A)\dlto^{j_!}\\
&\bK_G(TV;A)&
\enddiagram
 $$
or
 $$
\diagram
\bK_G(TX;A)\rrto^{i_!}\drto^{(j\circ i)_!\:}\ddrto_{\ti^X_G\:}
&&\bK_G(TY;A)\dlto_{j_!}\ddlto^{\,\ti^Y_G}\\
&\bK_G(TV;A)&\\
&K^G(A).\uto_\f&
\enddiagram
 $$
\hfill$\Box$

\dfi
An
{\em index function\/}
\label{ind:indexfunction}
is a family of $R(G)$-ho\-mo\-mor\-phisms
$\{\ind^X_G\}$
 $$
\ind^X_G:\bK_G(TX;A)\to K^G(A),
 $$
where $G$ runs through the set of compact Lie groups and $X$ is a smooth
compact $G$-manifold.
This family is restricted to satisfy
the following two conditions:
\begin{enumerate}
  \item If $f: X\to Y$ is a $G$-diffeomorphism, then
	the diagram
 $$
\diagram
\bK_G(TX;A)\rrto^{(df^{-1})^*}\drto_{\ind^X_G\:}
&&\bK_G(TY;A)\dlto^{\ind^Y_G}\\
&K^G(A)&
\enddiagram
  $$
is commutative.

  \item If $\psi: H\to G$ is a homomorphism of groups, then
	the diagram
 $$
\diagram
\bK_G(TX;A)\rto^{\psi^*}\dto_{\ind^X_G}&\bK_H(TX;A)\dto^{\ind^X_H}\\
\bK^G(A)\rto^{\psi^*}&\bK^H(A)
\enddiagram
 $$
is commutative.
 \end{enumerate}

 \begin{ass}
The topological index $\ti^X_G$ is an index function.
 \end{ass}

\Dok
1). Suppose $i:Y\hookrightarrow V,\: j=i\circ f:X\hookrightarrow V$.
By \ref{th:prophys}.2, the following diagram is commutative
 $$
\diagram
\bK_G(TY;A)\drto^{i_!}&&\\
&\bK_G(TV;A)\rto^{\:\f^{-1}}&K^G(A).\\
\bK_G(TX;A)\uuto^{f_!}\urto_{j_!}&&
\enddiagram
 $$
By \ref{th:prophys}.4, we have in this case $f_!=(df^{-1})^*$ and it
remains to use the definition of $\ti$.

2). Immediately follows from the definition. \hfill $\Box$

\smallskip
Let us consider the following two axioms.

{\bf Axiom A1.\/}\label{a1}
If $X={\rm one point}$, then
$\ind^X_G:\bK_G(TX;A)\to K^G(A)$
coincides with $\Id_ {K^G (A)} $.

{\bf Axiom A2.\/}\label{a2} Suppose $i: X\to Y$ is a $G$-enclosure, then
the diagram
 $$
 \diagram
\bK_G(TX;A)\rrto^{i_!}\drto_{\ind^X_G\:}&&\bK_G(TY;A)\dlto^{\ind^Y_G}\\
&K^G(A)&\\
\enddiagram
 $$
is commutative.

 \begin{cor}
{\em (from Theorem~\ref{th:ax2})\/} The topological index $\ti^X_G$
satisfies axioms {\em A1\/} and {\em A2.\/}
\hfill$\Box$
 \end{cor}

 \begin{th}\label{th:ax6}
Let $\ind^X_G$ be an index function satisfying axioms
{\em A1\/} and {\em A2\/}. Then $\ind^X_G=\ti^X_G$.
 \end{th}

\Dok
Consider a $G$-enclosure $i: X\to V$ of a manifold $X$ in a real
vector $G$-space $V$. The one-point compactification $V^+$
(i.~e., a sphere) is a $G$-manifold with the canonical $G$-inclusion
$\e^+:V\to V^+$. Suppose $i^+=\e^+\circ i:X\to V^+$.
If $P=0\in V$ and $j: P\to V$ is the inclusion, then
we obtain the diagram
 $$
 \diagram
&\bK_G(TX;A)\dlto_{i_!}\drto^{\ind^X_G\:}\dto^{i^+_!}&\\
\bK_G(TV;A)\rto^{(\e^+)_!}&\bK_G(TV^+;A)\rto^{\ind^{V^+}_G}&K^G(A)\\
&\qquad \bK_G(TP;A)=K^G(A),\ulto^{j_!}\uto_{j^+_!}\urto_{\ind^P_G}
\enddiagram
 $$
where $j^+=\e^+\circ j:P\to V^+$.
By \ref{th:prophys}.2 (resp., by axiom A2), the left (resp., right)
triangles commute.
By A1, $\ind^P_G$ is the identity mapping. Hence
\begin{eqnarray*}
\ind^X_G &=& \ind^{V^+}_G\,i^+_!= \ind^{V^+}_G\,(\e^+)_!\,i_!=
\ind^{V^+}_G\,(j^+)_!\,j^{-1}_!\,i_!=\\
&=& \ind^P_G\,j^{-1}_!\,i_!= j^{-1}_!\,i_!=\ti^X_G.
\end{eqnarray*}
Since
 $
j_!:K^G(A)\to \bK_G(TV;A)=\bK_G(V\o\C;A)
 $
coincides with the Thom homomorphism, the theorem is proved.
\hfill$\Box$

\smallskip
{\bf Axiom B1 (excision).\/} Let $U$ be a (noncompact)
$G$-manifold and
 $$
j_1:U\to X_1,\qquad j_2:U\to X_2
 $$
be $G$-enclosures of the manifold $U$ on open subsets of compact
$G$-manifolds $X_1$ and $X_2$. Then the diagram
 $$
 \diagram
&\bK_G(TX_1;A)\drto^{\ind^{X_1}_G}& \\
\bK_G(TU;A)\urto^{(dj_1)_*}\drto_{(dj_2)_*}&&K^G(A)\\
&\bK_G(TX_2;A)\urto_{\ind^{X_2}_G}&
\enddiagram
 $$
is commutative.

\smallskip
Suppose there exists though one of
indicated enclosures. Then by the axiom, the index
 $$
\ind^U_G:\bK_G(TU;A)\to K^G(A)
 $$
can be well defined.

Let us denote by $\ci^Y_H$ the classical (complex) index
 $$
\ci^Y_H:K_H(TY)\to R(H).
 $$
We have the following statement (see \cite{AZ1}).
Suppose $j:*\to\R^n$ is the enclosure of $\vec 0$, hence
 $j_!:R(O(n))\to K_{O(n)}(T\R^n)$. Then
$\ci^{\R^n}_{O(n)}\,j_!(1)=1$.

Let $\pi: P\to X$ be a compact differentiable principal
bundle for a group $H$ (compact Lie group).
Then we have (right) free action of $H$ on $P$ and $X=P/H$.
Suppose we have left action of $G$ on $P$ and these two
actions commute. Let $F$ be a compact left
$(G\times H)$-manifold. We can form the associated bundle
$\pi_1: Y=P\times_H F\to X$ with the natural action of $G$.
Consider the tangent bundle along the fibers of $\pi_1$.
Let us denote it by $T_F Y$. Then $T_F Y$ is a
$G$-invariant real subbundle of $TY$
and $T_F Y=P\times_H TF$. Using the metric it is
possible to decompose $TY$ into a direct sum
 $TY=T_F Y\oplus \pi_1^*(TX)$.
Therefore the multiplication
 $$
 \begin{array}{ccc}
 \bK_G(TX;A)\o K_G(T_F Y)&&\\ \downarrow &&\\
 \bK_G(\pi_1^*TX;A)\o K_G(T_F Y)&\to &\bK_G(TY;A)
 \end{array}
 $$
is defined.
There exists the map
 $$
K_{G\times H}(TF)\to K_{G\times H}(P\times TF)\cong
K_{G}(P\times_H TF)=K_{G}(T_F Y).
 $$
Hence we can define the mapping
 $$
\g:\bK_G(TX;A)\o K_{G\times H}(TF)\to \bK_G(TY;A).
 $$
Let us denote
$\g (a\o b)$ by $a\cdot b$.

If $V$ is a complex vector $(G\times H)$-space, then
$P\times_H V$ is a complex vector $G$-bundle over $X$.
We obtain the following ring homomorphism being
a homomorphism of $R(G)$-modules:
 $$
\mu_P:R(G\times H)\to K_G(X),\qquad [V]\mapsto [P\times_H V].
 $$
Since $\bK_G (TX; A) $ has a
$K_G (X)$-module structure, we can formulate the following axiom.

{\bf Axiom B2.\/}\label{b2}
If $a\in \bK_G(TX;A),\: b\in K_{G\times H}(TF)$, then
 $$
\ind^Y_G(a\cdot b)=\ind^X_G(a\cdot \mu_P(\ci^F_{G\times H}(b))),
 $$
i.~e., the diagram
 $$
 \diagram
\bK_G(TX;A)\otimes K_{G\times H}(TF)\qquad
\dto^{\g}\rto^{1\,\o\,\ci^F_{G\times H}}
& \qquad \bK_G(TX;A)\otimes R(G\times H)\dto^{1\o\mu_P}\\
\bK_G(TY;A)\dto^{\ind^Y_G}&\bK_G(TX;A)\otimes K_{G}(X)\dto \\
K^G(A)&\bK_G(TX;A)\lto_{\ind^X_G}
\enddiagram
 $$
is commutative

 \begin{th}\label{th:ax7}
Let $\pi: P\to X$ be a principal right $H$-bundle with a left action
of $G$ commuting with $H$. Suppose $F$ is a $(G\times H)$-space.
Let us denote by $Y$ the space $P\times_H F$. Let $j: X_1\to X$ and
$k: F_1\to F$ be $G$- and $ (G\times H)$-enclosures, respectively;
let $\pi^1: P_1\to X_1$ be the principal
$H$-bundle induced by $j$ on $X_1$; assume $Y_1:=P_1\times_H F_1$.
The enclosures $j$ and $k$ induce
$G$-enclosure $j*k: Y_1\to Y$. In this situation the diagram
 $$
 \diagram
\bK_G(TX;A)\otimes_{R(G)} K_{G\times H}(TF)\rto^{\qquad\qquad \g} &
\bK_G(TY;A)\\
\bK_G(TX_1;A)\otimes_{R(G)} K_{G\times H}(TF_1)\rto^{\qquad \qquad \g }
\uto_{j_!\o k_!}
& \bK_G(TY_1;A)\uto_{(j*k)_!}
\enddiagram
 $$
is commutative.
 \end{th}

\Dok
Let us describe $\g$ explicitly:
 $$
 \begin{array}{cc}
\bK_G(TX;A)\otimes K_{G\times H}(TF)\to
&\hspace{-2em} \bK_G(TX;A)\o K_{G\times H}(P\times TF)\cong\\
\qquad \uparrow\,j_!\o
 k_!\qquad \qquad\fbox{1}&\qquad\uparrow\,\e\qquad\qquad\fbox{2}\\
\bK_G(TX_1;A)\otimes K_{G\times H}(TF_1)\to
&\hspace{-1em} \bK_G(TX_1;A)\o K_{G\times H}(P_1\times TF_1)\cong
 \end{array}
 $$
\smallskip
$$
 \begin{array}{c}
\cong \bK_G(TX;A)\o K_G(P\times_H TF)\to\\
\fbox{2}\qquad \b\,\uparrow\qquad \qquad \fbox{3}\\
\cong \bK_G(TX_1;A)\o K_G(P_1\times_H TF_1)\to
 \end{array}
$$
\smallskip
\begin{equation}\label{axio1}
 \begin{array}{c}
\to \bK_G(\pi^*_1 TX;A)\o K_G(P\times_H TF)\to\\
\fbox{3}\qquad \qquad \a\,\uparrow \qquad \qquad \fbox{4}\\
\to \bK_G((\pi^1)^*_1 TX_1;A)\o K_G(P_1\times_H TF_1)\to
 \end{array}
\end{equation}
\smallskip
 $$
 \begin{array}{c}
\to \bK_G((\pi^*_1  TX)\times(P\times_H TF);A) \to \\
	     \fbox{4} \\
\to \bK_G((\pi^1)^*_1  TX_1\times(P_1\times_H TF_1);A) \to
 \end{array}
 $$
\smallskip
 $$
 \begin{array}{c}
\to \bK_G(\pi^*_1 TX\oplus (P\times_H TF);A)=\bK_G(TY;A)\\
\qquad \qquad \fbox{4}\qquad \qquad\qquad \qquad
\qquad \qquad \uparrow\, (j*k)_! \\
\to \bK_G((\pi^1)^*_1 TX_1\oplus (P_1\times_H TF_1);A)=\bK_G(TY_1;A).
 \end{array}
 $$
Let us remind the diagram, which was used
 for the definition of the Gysin homomorphism of an enclosure
$j: X_1\to X$:
 $$
\diagram
p^*_T(N_{X_1}\oplus N_{X_1})\dto &&TN_{X_1}\llto_\psi\rto^{d\F_{X_1}}
\dto& TW_{X_1} \rto^{dj_2}\ddto& TX\ddto \\
TX_1 \drto^{p_T} & & N_{X_1}\dlto_{p_{N_{X_1}}\:} \drto^{\F_{X_1}}&& \\
 & X_1\rrto^{j_1} & & W_{X_1} \rto^{j_2} & X.
\enddiagram
 $$
From the similar diagrams for $k_!$ and $ (j * k) _!$
and the explicit form of the maps it follows that
the square \fbox{4} in (\ref {axio1}) iff
$\a$ has the following form:
 $$
 \begin{array}{l}
\a(\s\o \rho)=(\pi^*_1)\Bigl\{
(dj_2)_*\,(d\F^{-1}_{X_1})^*\,\psi^*_{X_1}\Bigr\}\circ\f^S_A(\s)\o\\
\qquad \o(\pi^*j_2\times_H  dk_2)_*\,\Bigl((\pi^*\F_{X_1}\times_H
d\F_{F_1})^{-1}\Bigr)^*\,(1\times_H\psi_{F_1})^*\,\f^T_\C(\rho),
 \end{array}
 $$
where $S$ and $T$ are bundles of the form
 $$
 \begin{array}{cccc}
 &\pi^*_1\,\Bigl((p^X_T)^*\{N_{X_1}\oplus N_{X_1}\}&&
\pi^* N_{X_1}\times_H (p^{F_1}_T)^*\,(N_{F_1}\oplus N_{F_1})\\
S:&\qquad \downarrow\,(\pi^1)^*_1q^{X_1}_T &
T:&\qquad \downarrow\,(\pi^1)^*(p_{N_{X_1}})\times_H q^{F_1}_T \\
&(\pi^1)^*_1\,(TX_1), &&
\pi^*\,X_1\times_H TF_1=P_1\times_H TF_1.
 \end{array}
 $$
Hence the square \fbox{3} in (\ref {axio1}) is commutative iff the
homomorphism $\b$ has the form
 $$
 \begin{array}{l}
\b(\t\o\rho)=j_!(\t)\o \\
\:\o(\pi^*j_2\times_H dk_2)_*
\Bigl((\pi^*\F_{X_1}\times_H
d\F_{F_1})^{-1}\Bigr)^*\,(1\times_H\psi_{F_1})^*\,\f^T_\C(\rho),
 \end{array}
  $$
 $$
\t\in \bK_G(TX_1;A),\qquad \rho\in K_G(P_1\times_H TF_1).
 $$
In turn, the square \fbox{2} in (\ref {axio1})
is commutative iff the homomorphism $\e$ has the form
 $$
 \begin{array}{l}
\e(\t\o\delta)=j_!(\t)\o \\
\:\o(\pi^*j_2\times_H dk_2)_*
\Bigl((\pi^*\F_{X_1}\times_H
d\F_{F_1})^{-1}\Bigr)^*\,(1\times_H\psi_{F_1})^*\,\f^{\widetilde
T}_\C(\delta),
\end{array}
  $$
 $$
\t\in \bK_G(TX_1;A),\qquad \delta\in K_{G\times H}(P_1\times_H TF_1),
 $$
where $\widetilde T$ is the following bundle:
 $$
 \begin{array}{cc}
&\pi^* N_{X_1}\times (p^{F_1}_T)^*\,(N_{F_1}\oplus N_{F_1})\\
\widetilde T:\quad&\qquad \downarrow\,(\pi^1)^*(p_{N})\times
q^{F_1}_T \\
&P_1\times TF_1.\qquad \qquad
 \end{array}
 $$
Suppose $\delta=[{\underline \C}]\widehat\o\omega$, where
$[{\underline \C}]\in K_{G\times H}(P_1)$, ${\underline \C}$ is
the one-dimensional trivial bundle and $\omega\in K_{G\times H}(TF_1)$.
Then
\begin{eqnarray*}
\e(\t\o\delta)&=&j_!(\t)\o\Bigl\{\pi^*(j_2)_*(\F^{-1}_{X_1})^*
[{\underline \C}]\widehat\o k_!(\omega)\Bigr\}=\\
&=& j_!(\t)\o\Bigl\{[{\underline \C}]\widehat\o k_!(\omega)\Bigr\}.
\end{eqnarray*}
Since the map $K_{G\times H}(TF)\to K_{G\times H}(P\times TF)$ (as well
as the lower line in (\ref{axio1})) has the form
 $
\omega\mapsto [{\underline \C}]\widehat\o\omega ,
 $
we have proved the commutativity of \fbox {1} in (\ref{axio1}).
 \hfill$\Box$

\smallskip
Let in the situation of axiom B2
 $$
\ci^F_{G\times H}(G)\in R(G)\ss R(G\times H).
 $$
Since $\mu_P$ and $\ind^X_G$ are $R(G)$-homomorphisms,
the following property holds.

\smallskip
{\bf Axiom B2${}'$.\/} (corollary of B2) \ If
$\ci^F_{G\times H}(G)\in R(G)\ss R(G\times H)$, then
 $$
\ind^Y_G(a\cdot b)=\ind^X_G(a)\cdot \ci^F_G(b).
 $$

Assume in B2 $X=P,\: H=1$. We obtain the following axiom.

\smallskip
{\bf Axiom B2$''$.\/} (corollary of B2) \ If $X$ and $F$ are
$G$-manifolds, then
 $$
\ind^{X\times F}_G(a\cdot b)=\ind^X_G(a)\cdot \ci^F_G(b).
 $$

 \begin{th}\label{ax8}
Suppose an index function $\ind^X_G$ satisfies
{\em A1, B1, B2,\/} then
 $$
\ind^X_G=\ti^X_G.
 $$
 \end{th}

\Dok
Suppose in the axiom B2$'$  $F$ is equal to an open
$(G\times H)$-subset of the compact manifold $\widetilde F$. Let
$j: F\hookrightarrow\widetilde F$. Then
\begin{eqnarray}
\ind^Y_G(a\cdot b)&=&\ind^{\widetilde Y}_G\,(dJ_*)(a\cdot b)=
\ind^{\widetilde Y}_G\,(a\cdot ((dj)_* b))=\nonumber\\
&=&\ind^X_G(a)\cdot \ci^{\widetilde F}_{G\times H}((dj)_* b)=
\nonumber\\
&=&\ind^X_G(a)\cdot \ci^{F}_{G\times H}(b),\label{axio2}
\end{eqnarray}
where $J$ is the enclosure
 $$
Y=P\times_H F\:\stackrel{\Id\times_H j}{\hookrightarrow}\:P\times_H
\widetilde F=\widetilde Y.
 $$
Indeed, let us consider the diagram
 $$
 \begin{array}{c}
\bK_G(TX;A)\otimes K_{G\times H}(TF)\to
\bK_G(TX;A)\o K_{G\times H}(P\times TF)\cong\\
\qquad \downarrow\,1\o (dj)_*
 \qquad \qquad\qquad\downarrow\,1\o (\Id\times dj)_*\\
\bK_G(TX;A)\otimes K_{G\times H}(T\widetilde F)\to
\bK_G(TX;A)\o K_{G\times H}(P\times T\widetilde F)\cong
 \end{array}
 $$
\smallskip
 $$
 \begin{array}{c}
\cong
\bK_G(TX;A)\o K_{G}(P\times_H TF)=
\bK_G(TX;A)\o K_{G}(T_F Y)\to\\
\qquad \downarrow\,1\o (\Id\times_H dj)_*
 \qquad \qquad\qquad\downarrow\,1\o\a_* \\
\cong
\bK_G(TX;A)\o K_{G}(P\times_H T\widetilde F)=
\bK_G(TX;A)\o K_{G}(T_{\widetilde F} Y)\to
 \end{array}
 $$
\smallskip
 $$
 \begin{array}{c}
\to
\bK_G(\pi^*_1 TX;A)\o K_{G}\o K_{G}(T_F Y)\to \bK_G(TY;A)\\
\qquad \qquad \downarrow\,1\o \a_*
\qquad  \qquad \qquad\qquad\downarrow\,(dJ)_* \\
\to
\bK_G(\pi^*_1 TX;A)\o K_{G}(T_{\widetilde F} Y)\to
\bK_G(T\widetilde Y;A).
 \end{array}
 $$
This diagram is commutative. In fact, we have
 $$
 \begin{array}{ccc}
 TY&=&T_F Y\oplus \pi^*_1(TX)\\
\quad\downarrow\,dJ& &\qquad \downarrow\,{\scriptstyle \Mat \a 001} \\
 T\widetilde Y&=& T_{\widetilde F}Y\oplus \pi^*_1(TX)
 \end{array}
 \; ,
 $$
and $\a=\Id\times_H dj$ under the identification $T_F Y=P\times_H TF$.
We have proved the second equality in (\ref{axio2}), the remaining are
obvious.

Let us now take, in particular,
 $$
F=\R^n,\quad \widetilde F=(\R^n)^+=S^n,\quad H=O(n), \quad b=j_!(1),
\quad 1=[\underline\C],
$$
where $j: \vec 0\hookrightarrow \R^n$
is the natural enclosure. Then $P$ is a principal
$O(n)$-bundle over $X$, the group $G$ acts on $P$ commuting with
$O(n)$. Suppose $G$ acts on $\R^n$ in a trivial way.
We form the associated real $G$-bundle
 $$
P\times_{O(n)} \R^n=Y\to X.
 $$
Let us denote by
 $$
i:X\to Y,\qquad i=1_X * j,
 $$
the enclosure of $X$ as of the zero section. Assume
  in Theorem~\ref{th:ax7} $F=\R^n$, $X_1=X$, $F_1=pt$. Then
we obtain the commutative diagram
 $$
 \begin{array}{ccc}
\bK_G(TX;A)\o_{R(G)} K_{G\times O(n)}(T\R^n)&\ola\g& \bK_G(TY;A)\\
\qquad \uparrow\,(1_X)_!\o j_!&
& \uparrow\, (1_X * j)_!\lefteqn{=i_! }\\
\bK_G(TX;A)\o_{R(G)} K_{G\times O(n)}(pt)&\ola\g& \bK_G(TX;A).
 \end{array}
 $$
Since $\g\,(a\o 1)=a$,
\begin{eqnarray*}
i_!(a)&=&\g\,\Bigl(((1_X)_!\o j_!)\,(a\o 1)\Bigr)=\g\,(a\o j_!(1))=\\
&=&a\cdot j_!(1)=a\cdot b.
\end{eqnarray*}
By the property of $\ci$ indicated above,
 $$
\ci^{\R^n}_{O(n)}\,j_!(1)=1,\qquad \ci^{\R^n}_{G\times O(n)}\,j_!(1)=1,
 $$
where $G$ acts on $\R^n$ in a trivial way. Now by the axiom B2,
\begin{eqnarray}
\ind^X_G&=&\ind^X_G(a\cdot 1)=\ind^X_G(a\cdot\mu_P(1))=\nonumber\\
&=&\ind^X_G(a\cdot\mu_P(\ci^{\R^n}_{G\times O(n)}(b)))=\nonumber\\
&=&\ind^Y_G(a\cdot b)=\ind^Y_G \,i_!(a).\label{axio3}
\end{eqnarray}

Let $k: X\to Z$ be an enclosure of $X$ in a compact $G$-manifold
$Z$ with the normal bundle $N$ and a tubular neighborhood
$\F: N\to W$. By the definition of the Gysin homomorphism,
$k_!=(di_2\circ d\F)_*j_!$, where
$di_2: TW\to TZ$ is an enclosure of tangent bundle and $j: X\to N$
is the enclosure of $X$ as of the zero section in the normal bundle.
In the diagram
 $$
\diagram
\bK_G(TX;A)\rto^{j_!}\drto_{\ind^X_G\:}&
\bK_G(TN;A)\quad\rto^{(di_2\circ d\F)_*} \dto^{\ind^N_G}&
\quad \bK_G(TZ;A)\dlto^{\ind^Z_G}\\
&K^G(A)&
\enddiagram
 $$
the left triangle is commutative by (\ref {axio3}).
The map $i_2\cdot\F$ is an open enclosure. Hence by B2, the right
triangle is commutative too. Therefore,
 $
\ind^X_G=\ind^Z_G\circ k_!
 $
Hence A2 is satisfied. To complete the proof it remains to apply
Theorem~\ref{th:ax6}.
\hfill$\Box$

 \section{Proof of the index theorem}

First of all obviously, the analytical index is an index
function.

 \begin{lem}
The analytical index $\ai$ satisfies the axiom {\rm A1.}
 \end{lem}

\Dok
Elliptic operator over a point is a $GGA$-mapping
$P: V\to W$ of projective $GGA$-modules and
$[\s(P)]=[V]-[W]=\index P\in K^G(A).$\hfill$ \Box$

 \begin{th}
The index $\ai$ satisfies the axiom {\rm B1.}
 \end{th}

\Dok
Suppose  $a\in \bK_G(TU;A)$;
$j_1: U\to X_1$ and  $j_2: U\to X_2$
are $G$-enclosures; $\pi: TU\to U$ is the natural projection.
Let the sequence
 $$
0\to\pi^* E\ola\rho \pi^* F\to 0
 $$
be exact for $x\in U\setminus L$, $|\xi|>c$ (point and (co)vector).
Suppose
 $$
\a: E|_{U\setminus L}\cong (U\setminus L)\times N,
\qquad
\b: F|_{U\setminus L}\cong (U\setminus L)\times N,
 $$
 $$
\rho=(\pi^*\b)^{-1}\,(\pi^*\a),
 $$
$L$ is some $G$-invariant compact set.
Then it is possible to assume
the symbols
 $
\s_1\in\Smbl_0(X_1,E_1,F_1),\; \s_2\in\Smbl_0(X_2,E_2,F_2)
 $
be as follows. Suppose
 $$
E_1=E\cup_{j_1\a}(X_1\setminus j_1L)\times N,\quad
E_2=E\cup_{j_2\a}(X_1\setminus j_2L)\times N.
 $$
Let
the similar equalities hold for $F_1$ and $F_2$, and
 $\s_1=\rho\cup_{j_1}\Id$, $\s_2=\rho\cup_{j_2}\Id.$
Let us pass to construction of operators
$\widetilde D_1$ and $\widetilde D_2$,
which represent these symbols
in $\Int_0(X_1;E_1,F_1)$ и $\Int_0(X_2;E_2,F_2)$, respectively.
Let us take a trivializing cover, a partition of unity and
smoothing functions on $U$.
Pull them back on $j_1 U$ and $j_2 U$, and then complete these
collections of open sets (to obtain covers)
by some open sets not intersecting with
$j_1 L$ and $j_2 L$, respectively.
By our symbols and with the help of this data
let us construct in the usual way (non-invariant)
operators $D_1,\, D_2\in \CZ_0$, and then
 $$
\widetilde D_1={\rm Av}_G D_1\in \Int_0(X_1),\quad
\widetilde D_2={\rm Av}_G D_2\in \Int_0(X_2).
 $$
It is necessary to check up the equality
 $$
\index \widetilde D_1=\index \widetilde D_2\in K^G(A).
 $$
Since $L$ is invariant, the averaging over
this set is the same for both operators.
Since the operators have the order 0, we compute index in
$L_2$-spaces. For these spaces
 $$
L_2(X_1,E_1)\cong L_2(j_1 L, E_1|_{j_1 L})\oplus
L_2(X_1\setminus j_1 L, E_1|_{X_1\setminus j_1 L})
 $$
and
 $$
\widetilde D_1:L_2(X_1\setminus j_1 L, E_1|_{X_1\setminus j_1 L})\cong
L_2(X_1\setminus j_1 L, E_1|_{X_1\setminus j_1 L})
 $$
(this is the identity operator). Similar relations hold for
$\widetilde D_2$.  On the second summand of the decomposition of $L_2$
we have the commutative diagram
 $$
\diagram
\G(E_1|_{j_1 L})\rto^{\widetilde D_1}\dto_{(j_2j_1^{-1})}^\cong &
\G(F_1|_{j_1 L})\dto^{(j_2j_1^{-1})}_\cong \\
\G(E_2|_{j_2 L})\rto^{\widetilde D_2}&
\G(F_2|_{j_2 L}).
\enddiagram
 $$
This diagram demonstrates the coincidence of indices.
\hfill$\Box$

 \begin{th}
The analytical index $\ai$ satisfies the axiom {\rm B2}.
 \end{th}

\Dok
Consider the manifold $Y=P\times_H F$ over a compact
manifold $X$, where $P\to X$ is the principal bundle with
compact Lie group $H$, and $F$ is a compact left $H$-manifold.
The compact Lie group $G$ acts on $X=P/H$ and $Y=P\times_H F$.
Hereinafter the metrics on $P$, $X$, $F$, $Y$ and on
$A$-vector bundles are supposed to be invariant.
Let us recall that we considered in the axiom B2
two elements and their product of special form:
 $$
a\in \bK_G(TX;A),\quad b\in K_{G\times H}(TF),\quad a\cdot b\in
\bK_G(TY;A) .
 $$

Let us represent these elements by symbols of elliptic operators.
Let $a$ be represented by a smooth $G$-symbol of order 1.
Let $\A_1\in\CZ_1$ (not necessary invariant)
have the symbol $\a$.
Let us choose a trivializing cover $\{U_j\}$ of the manifold $X$
for $P$ and $Y$ with a subordinate partition of unity $\{\f_j^2\}$.
Let us consider the operator $\A^j_1=\f_j\A_1\f_j$ on $U_j$.
Suppose $p: Y\to X$ is the projection,
and $Y_j=p^{-1}(U_j)\cong U_j\times F$. The lifting
$\widetilde\A^j_1$ of $\A^j_1$ to $Y_j$ is a tensor
product. Hence it belongs to $\Int_1 (Y_j) $. Extend this lifting
by zero and get the element of $\Int_1 (Y) $. Let us average over $G$:
 $$
\widetilde\A={\rm Av}_G(\sum \widetilde\A^j_1)\in\Int_1(Y).
 $$
Then $\s_1 (\widetilde\A) =\widetilde\a$, where the lifting
$\widetilde\a$ of the symbol $\a$ is defined globally. This follows
from the invariance of the metric: for the decomposition of
cotangent space $Y$ into the vertical component $\eta$
horizontal component $\xi$ we assume
$\widetilde\a(\xi,\eta)=\a(\xi)$.

Let us restrict
$\wA^j_1$ on sections, which are constants along the fibers of $Y$,
then we receive $\A^j_1$. Therefore the restriction $\A$ on the space of
these sections is a $G$-invariant operator
$\A={\rm Av}_G(\sum\A^j_1)\in\CZ_1(X)$ with symbol $\a$.

Let
 $
b=[\b]\in K_{G\times H}(TF),\;\b=\s(\B),\; \B\in\CZ_1(F)
 $
is a $(G\times H)$-operator. Let $\wB_1$ be the operator over
$P\times F$, obtained by a lifting of $\B$. Since this operator is
$(G\times H)$-invariant, it induces a $G$-invariant operator
$\wB$ over $Y=P\times_H F$ by the restriction of $\wB_1$ on the constant
sections along fibers of $P\times F\to Y$.
Since $P$ is locally trivial, the restriction $\wB_j=\wB|_{U_j}$
over $Y_j=p^{-1}(U_j)\cong U_j\times F$ is the lifting of $\B$.
Then $\wB_j\in\Int_1(Y_j)$, $\wB\in\Int_1(Y)$. Suppose
 $
\widetilde\b=\s_1(\wB),\; \widetilde\b(\xi,\eta)=\b(\eta),
 $
where $\xi$ is the horizontal component and $\eta$ is the vertical one.
There is a $G$-invariant operator
 $$
D=\Mat \wA {-\wB^*}\wB{\wA^*}\in\Int_1(Y,
E^0\widehat\oplus G^0\,\oplus \, E^1\widehat\oplus G^1,
E^0\widehat\oplus G^1\,\oplus \,E^1\widehat\oplus G^0),
$$
$$
\s(D)=\Mat {\widetilde\a}{-\widetilde\b^*}{\widetilde\b}{\widetilde\a^*},
 $$
and $[\s (D)] = [\a] \,[\b] =a\cdot b$ by the definition of the
multiplication $\g$. Let us calculate now $\index D$. Since $\wB$ is
``of complex origins'', its kernel and cokernel are invariant modules
from $\P(A)$ and $\index\wB=[\Ker\wB]-[\Coker\wB]\in K^G(A)$. The
operator $\wB$ is the extension of $\B$ to fibers; $\Ker\wB$ consists
of those smooth sections, which lay in $\Ker\B_x$ for each fiber
$Y_x$.  Here $\B_x$ is the operator over $Y_x$ acting as $\B$ over
the standard fiber $F$.  Hence $\Ker\wB$ is the space of smooth
sections of the vector bundle $K_\B=P\times_H \Ker\B$ over $X$. Since
$\wA$ and $\wB$ commute, $\wA$  induces an operator $C^K$ on the
sections $K_\B$.  From the definition of $\wA$  it follows that
$$
C^K={\rm Av}_G\,\left(\sum \f_j C^K_j\f_j\right),
 $$
where $C^K_j$ is the operator induced by $\wA^j_1$ on $K_\B|_{U_j}$.
By the definition of $\wA^j_1$, this means that
$C^K_j=\A^j_1\o\Id_{K_\B}$.
Therefore
 $$
C^K_j\in\Int_1,\quad C^K\in\Int_1,\quad \s(C^K)=\a\o\Id_{K_\B}.
 $$
Thus $C^K$ is an elliptic $G$-invariant operator on $X$,
 $
[\s(C^K)]=a\,[K_\B]\in \bK_G(TX;A),
 $
and
 $$
\ai^X a[K_\B]=\index^X\,C^K\in K^G(A).
 $$
Similarly, if $L_\B=P\times_H\Coker\B$, then $\s(C^L)=\a\o\Id_{L_\B}$,
 $$
[\s(C^L)]=a\,[L_\B]\in \bK_G(TX;A),
 $$
 $$
\ai^X a[L_\B]=\index^X\,C^L\in K^G(A),
 $$
and
\begin{eqnarray*}
\index^X\,C^K-&&\hspace{-2em}\index^X\,C^L=\\
&&=\ai^X \Bigl(a([K_\B]-[L_\B])\Bigr)=\\
&&=\ai^X \Bigl(a\,\mu_P(\ci^F_{G\times H}\,b)\Bigr)\in K^G(A).
\end{eqnarray*}
It remains to show that
 $$
\index^Y \,D=\index^X\,C^K-\index^X\,C^L.
 $$
By the definition of the index of an $A$-Fredholm operator,
$$
\index^X\,C^K=[N^K_0]-[N^K_1],\quad C^K:M^K_0\cong M^K_1.
$$
Let us remark the following. If
 $$
F=\Mat {F_1}00{F_2},\quad F_1:M_0\cong M_1,\quad F_2:N_0\to N_1,
 $$
is a decomposition for $GGA$-Fredholm operator,
$\index F= [N_0] - [N_1] $, then $\index (F^*)=-\index F$ and
 $$
-\index^X\,C^L=\index^X\,(C^L)^*,
 $$
where $(C^L)^*$ is constructed in the same way as $C^L$, but instead of
the operator $\wA$ we take $\wA^*$. Suppose,
$$
\index^X(C^L)^*=[N^{L*}_0]-[N^{L*}_1],\quad (C^L)^*:M^{L*}_0\cong
M^{L*}_1.
  $$
Then by the definition of $K_\B$ and $L_\B$ as the kernels and cokernels,
and $C^K$ and $ (C^L)^*$, respectively, we have
 $$
 \begin{array}{cc}
\wA:M^K_0\cong M^K_1,&\wB(M^K_0)=0,\\
\wA^*:M^{L*}_0\cong M^{L*}_1,&\wB^*(M^{L*}_0)=0.
 \end{array}
 $$
Further,
\begin{equation}\label{exact1}
D\,\left(
 \begin{array}{c}
 M^K_0\\ \oplus \\ M^{L*}_0
 \end{array}
\right)
\,=\,
\left(
 \begin{array}{c}
\wA (M^K_0)\\ \oplus \\ \wA^*(M^{L*}_0)
 \end{array}
\right)
\,=\,
\left(
 \begin{array}{c}
 M^K_1\\ \oplus \\ M^{L*}_1
 \end{array}
\right),
\end{equation}
\begin{equation}\label{exact2}
D\,\left(
 \begin{array}{c}
 N^K_0\\ \oplus \\ N^{L*}_0
 \end{array}
\right)
\,=\,
\left(
 \begin{array}{c}
\wA (N^K_0)\\ \oplus \\ \wA^*(N^{L*}_0)
 \end{array}
\right)
\,\ss\,
\left(
 \begin{array}{c}
 N^K_1\\ \oplus \\ N^{L*}_1
 \end{array}
\right).
\end{equation}
This is the description of the action of
$D$ on $\Ker\wB\oplus \Ker\wB^*$.

Let now $x\in\Ker\wB^\bot$, i.~e.,
$(x,y)=0$ for any $y$ such that $\wB y=0$.
For these $x$ and $y$ we have
 $$
\wB(\wA^* y)=\wA^*(\wB y)=0,\quad (\wA x,y)=(x,\wA^* y)=0.
 $$
Therefore $\A\Bigl((\Ker\wB)^\bot\Bigr)\ss (\Ker\wB)^\bot$.
 Similarly,
$$
\A^*\Bigl((\Coker\wB)^\bot\Bigr)\ss (\Coker\wB)^\bot.
$$
Hence,
\begin{equation}\label{exact2a}
D\,\left(
 \begin{array}{c}
(\Ker\wB)^\bot\\ \oplus \\ (\Coker\wB)^\bot
 \end{array}
\right)
\,\ss\,
\left(
 \begin{array}{c}
(\Ker\wB)^\bot\\ \oplus \\ (\Coker\wB)^\bot
  \end{array}
\right).
\end{equation}
From
 $$
D^*D=\Mat {\wA^*\wA+\wB^*\wB}00{\wA\wA^*+\wB\wB^*} \ge
\Mat {\wB^*\wB}00{\wB\wB^*}
 $$
it follows that the operator $D$ is an isomorphism on
$(\Ker\wB)^\bot\oplus (\Coker\wB)^\bot$. Taking into account
the inclusion (\ref{exact2a}),
(\ref{exact1}) and (\ref{exact2}), we obtain the following
formula for the inverse image:
\begin{equation}\label{exact3}
D^{-1}\,\Bigl((\Ker\wB)^\bot\oplus (\Coker\wB)^\bot\Bigr)\,\ss\,
(\Ker\wB)^\bot\oplus (\Coker\wB)^\bot.
\end{equation}
To prove that
 $$
D:(\Ker\wB)^\bot\oplus (\Coker\wB)^\bot\,\to\,
(\Ker\wB)^\bot\oplus (\Coker\wB)^\bot
 $$
is an isomorphism, it remains to check up that it is an epimorphism.
By calculations as above, we obtain that $(\Ker\wB)^\bot$ is
an $\wA^*$-invariant submodule and $(\Coker\wB)^\bot$
is an $\wA$-invariant submodule.
The composition $DD^*$ defines the isomorphism
 $$
DD^*:(\Ker\wB)^\bot\oplus (\Coker\wB)^\bot\,\cong\,
(\Ker\wB)^\bot\oplus (\Coker\wB)^\bot.
 $$
By (\ref{exact3}), $D$ is an epimorphism.

Now we can calculate the index:
\begin{eqnarray*}
\index^Y\,D&=& [N^K_0\oplus N^{L*}_0]-[N^K_1\oplus N^{L*}_1]=\\
&=&([N^K_0]-[N^K_1])+([N^{L*}_0]-[N^{L*}_1])=\\
&=&\index^X\,C^K+\index^X\,C^{L*}=\\
&=&\index^X\,C^K-\index^X\,C^{L}=\\
&=&\ai^X\Bigl(a\cdot\mu_P(\ci^F_{G\times H}b)\Bigr).\qquad \Box
\end{eqnarray*}

 \begin{th}\label{Cind}
Index functions $\ai$ and $\ti$ coincide.
 \end{th}

\Dok
From the statements proved above, it follows that we can apply
Theorem \ref{ax8}.\hfill$\Box$

 \section{Families and algebras}\label{sec:famil}
 \markright{}

Let us remind the connection between
$C^*$-elliptic operators and families of elliptic operators for
non-equivariant case (cf.~\cite{MiUMN}). The equivariant case and
the case of families of $C^*$-operators can be obtained just similarly.

Let $M$ be a smooth closed manifold, $X$ be a compact Hausdorff space,
$A=C(X)$ is a unital $C^*$-algebra. Let us consider a family of
operators elliptic along the fibers of
$M\times X\to X$. More precisely, we have complex bundles
$E\to M\times X$, $F\to M\times X$ and an operator $D$, which is
continuous (e.~g. after enclosure of the bundles in a trivial one)
family of operators
$D_x:H^s(E_x)\to H^{s-m}(F_x)$. By Theorem of J\"anich \cite{Jani}
the index
 $$
\index_a \{D_x\}=[\Ker D_x] - [\Coker D_x]\in K^0(X)=K_0(A)
 $$
is defined (we get bundles after a compact perturbation).
On the other hand, $E$ can be considered as an $A$-bundle $E_A$ over $M$
(similarly, $F_M$), and we can define an operator
$D_A$ by commutativity of the following diagram
(we consider the continuous cross-sections of subbundle of product
bundle of $X$ by Sobolev space of stabilized fiber)
 $$
\diagram
\G(X,\{H^s(E_x)\})\rto^{\{D_x\}_*\ }&\G(X,\{H^{s-m}(F_x)\})\\
H^s(E_A)\rto^{D_A}\uto^{\p}_\cong&H^{s-m}(F_A)\uto^{\p}_\cong.
\enddiagram
 $$
So $\{D_x\}_*(h)(x):=D_x(h(x))$, where $h:X\to \{H^s(E_x)\}$, and
locally (when $E_A=E_x\ot A$)
$$
\p(f\ot a)(x)=f\cdot a(x),\qquad \p(D_A(f\ot a))(x)=D_x(f\cdot a(x))=
D_x(f)\cdot a(x).
$$
 \begin{lem}
The operator $D_A$ is an $A$-elliptic operator and
$$\ai D_A=\index_a \{D_x\}.$$
 \end{lem}

\Dok As it was demonstrated,
$D_A$ is an $A$-homomorphism. By calculation in local coordinates
we obtain immediately, that $D_A$ is an $A$-$\Psi$DO with the
correspondent symbol. The compactness of $X$ implies its invertibility
over the complement of a compact neighborhood in the cotangent bundle,
i.~e. ellipticity. The coincidence of (co)kernels follows from the
definition of $D_A$ after a compact perturbation, which has no
influence for indices.\qquad $\Box$

\end{document}